\documentclass[12pt]{article}
\usepackage{amsthm,amsmath,amsfonts,amssymb,array} 
\usepackage[latin1]{inputenc}

\newcommand{\be}{\begin{enumerate}}
\newcommand{\ee}{\end{enumerate}}
\newcommand{\beq}{\begin{eqnarray*}}
\newcommand{\eeq}{\end{eqnarray*}}
\newcommand{\beqnr}{\begin{eqnarray}}
\newcommand{\eeqnr}{\end{eqnarray}}

\newcommand{\beg}{\begin{equation}}
\newcommand{\ed}{\end{equation}}

\newcommand{\feddich}{\hfill $\Box$\\}

\newcommand{\R}{\mathbb{R}}
\newcommand{\N}{\mathbb{N}}

\newcommand{\supp}{\operatorname{supp}}

\newtheorem{assumption}{Property}
\newtheorem{proposition}{Proposition}

\newtheorem{theorem}{Theorem}
\newtheorem{lemma}{Lemma}

\newtheorem{remark}{Remark}
\newtheorem{corollary}{Corollary}

\def\T{{\relax\ifmmode I\!\!\hspace{-1pt}T\else$I\!\!\hspace{-1pt}T$\fi}}

\def\gsim{\mathrel{\raisebox{-4pt}{$\stackrel{\textstyle>}{\sim}$}}}

\def\lsim{\raisebox{-1ex}{$~\stackrel{\textstyle<}{\sim}~$}}

\def\cA{{\cal A}}

\def\bu{{\bf u}}

\def\bz{{\bf z}}

\def\cE{{\cal E}}

\def\cI{{\cal I}}

\def\cL{{\cal L}}
\def\cB{{\cal B}}

\def\cS{{\cal S}}

\def\bC{{\bf C}}

\def\bM{{\bf M}}

\def\bv{{\bf v}}
\def\bz{{\bf z}}
\def\bw{{\bf w}}

\def\lll{\langle}
\def\rr{\rangle}

\newcommand{\br}{{\bf r}}

\newcommand{\cX}{\mathcal{H}}
\newcommand{\XX}{\mathcal{X}}
\newcommand{\beqn}{\begin{equation}}
\newcommand{\eeqn}{\end{equation}}

\def\endproof{\hfill\rule{1.5mm}{1.5mm}\\[2mm]}

\newcommand{\bx}{\mathbf{x}}
\renewcommand{\bv}{\bx}
\newcommand{\by}{\mathbf{y}}
\newcommand{\bvv}{\mathbf{v}}

\makeatletter
\@addtoreset{equation}{section}
\makeatother

\newcommand\eref[1]{(\ref{#1})}

\def\int{\intop\limits}

\newcommand\bA{{\bf A}}
\newcommand\bB{{\bf B}}

\newcommand\bD{{\bf D}}

\newcommand\ve{\varepsilon}
\newcommand\cH{{\cal H}}

\newif\ifNZB

\makeatletter
\newcommand{\tr}{{\mathop{\operator@font T}\nolimits}}
\makeatother

\usepackage[usenames]{color}

\usepackage{algorithm}
\usepackage{algorithmic}

\begin{document}
\title{Adaptive Eigenvalue Computation - Complexity Estimates\thanks{
This research was supported in part by the Leibniz Programme of the DFG, by
the SFB 401 funded by DFG, the DFG Priority Program SPP1145 and by the
EU NEST project BigDFT.}
}
\author{W. Dahmen, T. Rohwedder, R. Schneider, A. Zeiser}
\date{November 7, 2007}
\maketitle

\begin{abstract}
This paper is concerned with the design and analysis of a fully
adaptive eigenvalue solver for linear symmetric  operators. After
transforming the original problem into an equivalent one formulated
on $\ell_2$, the space of square summable sequences, the problem
becomes sufficiently well conditioned so that a gradient type
iteration can be shown to reduce the error by some fixed factor per
step. It then remains to realize these (ideal) iterations within
suitable dynamically updated error tolerances. It is shown under
which circumstances the adaptive scheme exhibits in some sense
asymptotically optimal complexity.
\end{abstract}

\section{Introduction}\label{sect1}
\subsection{Background }\label{sect1.1}
In a Gelfand triple $\cH \overset{d}{\hookrightarrow} \XX
\overset{d}{\hookrightarrow}  \cH'$ of Hilbert spaces, where $ \cH $
is densely embedded in $ \XX$,   and  the dual pairing induced by
the inner product $\langle.,.\rangle$ of $\XX$, let
\begin{align}
  \cL: \cH \rightarrow \cH'
\end{align}
be a linear operator  that takes  $\cH$ onto its normed dual $\cH'$.
We wish to find, under certain assumptions on $\cL$,
 an eigenpair $(\lambda, u)$ of the problem
\beqn
\label{1.0}
\cL u=\lambda\cE  u
\eeqn
corresponding to the smallest eigenvalue $\lambda$.
Here $\cE: u \mapsto \langle \cdot, u \rangle$ canonically embeds $\cH$ in $\cH'$.
 Of course, \eref{1.0} is to be understood in the weak sense, i.e.
\beg
\label{1.2}
\lll v,\cL u\rr = \lambda \lll v, \cE u\rr,~~~~v\in \cX,
\ed
In view of the above choice of the dual pairing $\lll \cdot,\cdot\rr$ we shall
from now on simply write $\lll v,w\rr$ instead of $\lll v,\cE w\rr$ for $v,w\in \cH$.

For typical applications,
one can think of $\XX$ as an $L_2$-space over some domain $\Omega$ and $\cH$ as a Sobolev space $H^t$ of positive
real order $t$ (or as a closed subspace of a Sobolev space, determined e.g. by
homogeneous boundary conditions).

We shall be concerned with the case that
the operator $\cL$ is symmetric, i.e.
  \begin{align}
    \langle \cL v, w \rangle = \langle v, \cL w \rangle,
    \quad \textnormal{for all} \quad v, w \in \cH \label{Lsymmetric}
  \end{align}
and bounded and strongly positive, i.e. there exist positive constants $c_L, C_L$ such
that
\begin{align}
\label{Lelliptic}
c_L || v ||^2_\cH~~ \leq ~~\lll \cL v,v\rr =:|| v||_\cL^2 ~~\leq ~~C_L ||v||^2_\cH, \quad \quad v\in \cH.
\end{align}

Note that if the  (real) spectrum of $\cL$  is bounded from below by a negative number, $\cL$ can be shifted by a suitable parameter $\mu$ such that $\cL - \mu \cE$ satisfies \eref{Lelliptic}.  With this slight modification, this assuption is met by e.g. the electronic Schr\"odinger operator \cite{reed4}, one particle Schr\"odinger operators in $\R^3$ with certain potentials \cite{weidmann, hislob} and eigenvalue problems for strongly elliptic (satisfying a Garding inequality) differential operators on bounded domains as well as by eigenvalue problems for
Fredholm integral operators.

Note that \eref{Lelliptic} also  implies that $\cL$ is a norm-isomorphism from $\cH$ onto $\cH'$ whose condition is bounded by $C_L/c_L$.

 We shall exclusively deal with eigenvalue problems for which
 the infimum of the spectrum is an isolated eigenvalue $\underline{\lambda}$.
 For simplicity we shall actually assume that
$$
\underline \lambda := \inf \frac{\lll \cL u, u \rr}{\lll u,u \rr}
$$
is a simple eigenvalue with corresponding eigenvector $u$, and that that the rest of the spectrum is bounded from below by $\Lambda >\underline \lambda$, which means that
\beqnr
\label{specbound}
\frac{\langle \cL v, v \rangle}{\langle v, v \rangle} ~\geq ~\Lambda ~~~~\makebox {holds for all} ~~~ v \in \cH ~~\makebox {with} ~~~ \langle \cL u, v \rangle = 0. \eeqnr
The conventional approach is to
discretize \eref{1.0}, e.g. by finite elements or finite differences, which gives rise to
a finite dimensional discrete problem
\beqn
\label{1.0.1}
\bA_h \bu_h = \mu \bC_h \bu_h
\eeqn
where, for a given basis of the trial space, $\bA_h$, $\bC_h$ are the corresponding
stiffness matrix of $\cL$ and the mass matrix, respectively, and $\bu_h$ denotes the coefficient vector of the approximate eigenvector.
Now the issue becomes to solve \eref{1.0.1} efficiently within a suitable accuracy
tolerance associated with the discretization.
There is a vast literature on this issue, see e.g. Parlett \cite{parlett},
 Golub/Van Loan \cite{GvL}, which are  standard text books in numerical linear algebra,
or   Chatelin \cite{Chat},   Babuska-Osborn \cite{babo} concerning Galerkin discretizations and
     Dyaknov \cite{Dyak}, Knyazev/Neymair \cite{knyazev} for preconditioned iteration schemes. See also \cite{rsz}
for further references and comments.

In this paper we follow a different line. While the subsequent analysis will apply to
finite dimensional problems of the form \eref{1.0.1} as well, our primary interest
is to avoid the separation of the discretization and solution process. Rather, we will design
an abstract iteration scheme that solves the original {\em infinite dimensional} problem within a given accuracy tolerance, where the iterates of this scheme are taken, in principle, from all of $\cH$.
We will then show how to realize these updates (up to a perturbation) in a finite dimensional setting, particularly directing our attention to carrying out this task at possibly low computational cost.
 After completion of this work, we became aware of related work in \cite{GG}
 where an adaptive finite element scheme for elliptic eigenvalue problems is shown to converge without giving
 complexity estimates though.

In principle, such infinite dimensional iterations can be formulated
in different ways, focussing either on the convergence of Rayleigh
quotients or of eigendirections. While the first option - although
with the same motivation as in the present work - has been adopted
in \cite{rsz}, we address here the algorithmic realization of the
second option. In fact, in \cite{rsz}, from a somewhat different
perspective, we focus on  convergence of preconditioned  iteration
schemes per se,  whereas in the present paper we develop and analyze
a convergent  adaptive algorithm. Our analysis given here provides
(a posteriori) criteria for adaptive updates and yields complexity
estimates that prove the (asymptotic) optimality of the scheme.

We note that problems of the kind (\ref{1.2}) may be
 treated by means of  inverse iteration or better by Davidson or
Jacobi Davidson type methods as well. 
 These approaches require the solution of a
linear system in each iteration step. For this purpose
one may apply  
the adaptive solution strategies proposed in \cite{CDD1,CDD2,CDD3} based again on various types 
of (infinite dimensional) iteration schemes.
 In contrast to this strategy,
in the present approach we will intertwine the outer loop, e.g.
inverse iteration or Jacobi-Davidson and the iterative solver in the
inner loop by updating the Rayleigh quotient in each inner iteration step.
In fact, the present approach can be viewed as a preconditioned
steepest descent method for minimizing the Rayleigh quotient. 
At any rate, one should stress 
 that  for large  systems of linear
equations the use of iterative methods would be inevitable, anyway.

As in Jacobi Davidson or in Krylov space methods, one can also apply
subspace acceleration techniques to improve the present scheme by computing Ritz values of a
rather small system. This requires the computation of scalar
products like those appearing in the Rayleigh quotient.  Since we deal with infinite
matrices in the  present paper  we can only compute these values
approximatively, see also Section \ref{sect3.3} for further details.
Regarding these possible improvements, we content ourselves here with brief
indications and further studies are required.

We proceed now with showing how the problem \eref{1.0} is transformed into an
equivalent one formulated over a sequence space which is, however,
in some sense better conditioned. To this end, we begin with  collecting a few
facts that are relevant for carrying out this program. Moreover,
this will motivate the assumptions made in subsequent sections.

\subsection{Transformed Problems}\label{sect1.2}

Instead of projecting the problem \eref{1.0} to a (fixed) finite
dimensional subspace $\cH_h \subset \cH$ of $\cH$ ,
we shall transform
\eref{1.0} into an {\em equivalent} problem defined on the sequence
space $\ell_2(\cI)$ of square summable sequences of some possibly
infinite index set $\cI$ endowed with the inner product
$$
\langle \bvv, \bw \rangle = \sum_{i\in \cI} v_i w_i,\quad \textnormal{for all} \quad \bvv,\bw \in \ell_2(\cI)
$$
and induced norm $\|\cdot\|$. Note that we use the same symbols for the inner product on $\ell_2(\cI)$ as for the dual pairing on $\cH'\times\cH$.

The key ingredient is a suitable Riesz basis $\Psi =\{\psi_i:i\in \cI\}$ of $\cH$, i.e. there exist positive constants $c_\Psi, C_\Psi$ such that
\beqn
\label{ne}
c_\Psi ||\bvv||\leq ||\sum_{i\in \cI}v_i\psi_i||_\cH \leq C_\Psi ||\bvv||,
\quad \textnormal{for all} \quad \bvv \in \ell_2(\cI).
\eeqn
Then \eref{1.0} is equivalent to
\beqn
\label{1.4a}
\bA\bu = \lambda \bC \bu,  
\eeqn
where
\beqn
\label{1.4d}
 \bA:= \big( \lll \cL \psi_j,\psi_i\rr\big)_{i,j\in \cI}\,, \,\,\,\,\,
 \bC:= \big(\lll \psi_j,\psi_i\rr \big)_{i,j\in \cI}\,.
 \eeqn

For the operator $\bA$ the properties \eref{Lelliptic} of $\cL$ together with \eref{ne} imply (see e.g. \cite{Dacta}) that
\beqn
\label{1.4b}
||\bA||_{\ell_2 \to \ell_2} \leq C_L C_\Psi^2, \quad  \quad  ||\bA^{-1}||_{\ell_2 \to \ell_2}
\leq c_L^{-1} c_\Psi^{-2},
\eeqn
for the problem on $\ell_2(\cI)$, which in turn means that
\beqn
\label{1.4c}
c_L  c_\Psi^{2}||\bvv ||^2 \leq ||\bvv ||^2_\bA \leq C_L  C_\Psi^{2}||\bvv ||^2,
\eeqn
where $\|\cdot\|_{\bA}^2 := \langle \bA \cdot,\cdot\rangle$. Thus,
$\bA$-ellipticity is equivalent to bounded invertibility of $\bA$ on
$\ell_2(\cI)$. One could say that the transformation has a built-in
preconditioning effect: The original $\cL$ has been transformed to
an operator $\bA$ which is well conditioned on $\ell_2(\cI)$
according to \eref{1.4c}, a fact that will play a crucial role in
solving \eref{1.4a} numerically.

As for the properties of $\bC$, in the above setting
the operator $\bC$ stems from the inner product on $\XX$ and is therefore symmetric,  bounded and positive definite. The boundedness of $\bC$ as a linear operator on $\ell_2(\cI)$ results from the continuous embedding of $\cH$ in $\XX$.
However $\bC$ is typically not coercive on $\ell_2(\cI)$.

The role of $\bC$ is further illuminated when specifying the setting slightly
to the following situation which may actually serve as a guiding
example.
In fact, consider again the example $\cH =H^t, \XX =L_2$. Here, a suitable choice of the basis $\Psi$ is a wavelet-type basis.
In fact, such wavelet bases are available by now for a wide range of domains and exhibit an
important property, namely that a scaled version of $\Psi$ is also a Riesz basis for the pivot space
$\XX$. That means, setting
$$
\bD:= {\rm diag}(d_i:= \lll \psi_i,\psi_i\rr^{1/2}= ||\psi_i||_\XX : i\in \cI),
$$
 the collection
$$
\Psi^\circ:= \{d_i^{-1}\psi_i: i\in \cI\}
$$
satisfies
\beqn
\label{neX}
c_\Psi' ||\bD \bvv||\leq ||\sum_{i\in \cI}v_i\psi_i||_\XX \leq C_\Psi' ||\bD \bvv||, \quad \quad
\bD\bvv \in \ell_2(\cI),
\eeqn
where $c_\Psi', C_\Psi'$ are again fixed positive constants. It is well-known that for
$\XX =L_2$ and $\cH = H^t$ one has $d_i =\lll \psi_i,\psi_i\rr^{1/2} \sim 2^{-t|i|}$, where $|i|$
denotes the dyadic level of the wavelet $\psi_i$. Moreover,  the matrices
$\bA$ and $\bD^{-1}\bC\bD^{-1}$ are spectrally equivalent, i.e. there are constants
$c^*, C^*$ such that
\beqn
\label{c*}
c^*|| \bvv ||_{\bD^{-1}\bC\bD^{-1}}\leq ||\bvv||_\bA \leq C^*|| \bvv ||_{\bD^{-1}\bC\bD^{-1}}, \quad \,
\bvv \in \ell_2(\cI).
\eeqn
\begin{remark}
\label{remC}
In both of the above settings, $\bC$ is symmetric  positive definite and bounded but typically not coercive on $\ell_2(\cI)$.
However, $\bC$ is coercive on the space $\{\bv : \bD\bv\in \ell_2(\cI)\}$ equipped with the proper norm in the second example.

Note also that if the basis functions $\psi_i$ are pairwise orthogonal with respect to the
pivot inner product $\lll\cdot,\cdot\rr$, then $\bC$ is actually a diagonal matrix.
\end{remark}

The spectrum of the generalized eigenvalue problem $\bA \bu = \lambda \bC \bu$ coincides with the original spectrum of \eref{1.0}. In particular,
since $\bC$ is   positive definite (but not necessarily coercive) on $\ell_2(\cI)$, \ref{evp} is
equivalent to $\bC^{-1/2} \bA \bC^{-1/2} \by =\lambda \by$. The minimal eigenvalue
is then given by
\beqn \label{raylmin}
0 ~~<~~ \underline \lambda = \min_{v} \frac{\lll \cL v,v\rr}{ \lll v, v\rr}~~=~~
\min_{\by}\frac{\lll \bC^{-1/2} \bA \bC^{-1/2} \by,\by\rr}{ \lll \by,\by\rr}~~= ~~\min_{\bx}
\frac{\lll \bA\bx,\bx\rr}{\lll \bC\bx,\bx\rr},
\eeqn
so that the minimal eigenvalue of the transformed problem equals the one of the original problem
and $\underline \lambda$ is simple. Denoting by $\bu\in\ell_2(\cI)$ the corresponding eigenvector,
one has
\begin{align}
  \label{orthgap}
  \frac{\langle \bA \bvv,\bvv\rangle}{\langle \bC \bvv,\bvv\rangle} \ge \Lambda,
  \quad \textnormal{for all} \quad \bvv\in \ell_2(\cI) \quad \textnormal{with} \quad
  \langle \bA \bvv,\bu\rangle = 0.
\end{align}

 The properties collected above will guide us later when formulating the precise
 conditions all subsequent developments will be based upon.

 \subsection{Objectives and Layout}\label{sect1.3}

The above considerations indicate how to arrive at transformed
equivalent eigenproblems on sequence spaces \eref{1.4a} that are better
conditioned in the sense that the spectrum of the resulting matrix $\bA$ is enclosed in a bounded
interval of the positive real semi-axis. In particular, the
resulting ellipticity on certain orthogonal complements
will be seen to give rise to iteration schemes that reduce the error of the 
approximate smallest eigenvalue and also the deviation of the
corresponding approximate eigenspace from the exact one by at least
some fixed factor less than one. This is the case for the full
infinite dimensional problem as well as for any finite dimensional
discretization resulting from projecting onto spans of subsets of
the Riesz basis $\Psi$.

Once the fixed error reduction  has been established for the full infinite dimensional
problem, the principal idea is to mimic this ideal iteration on the infinite dimensional problem
(that still contains all information) by ultimately carrying out these iterations only
approximately. One then faces the following two tasks:
\begin{itemize}
\item[(i)] Find appropriate stage dependent tolerances within which each iteration is to be solved
so as to still guarantee convergence to the exact solution;
\item[(ii)]
 Devise numerical schemes
that realize the approximate application in the involved operators within the given
target accuracy at possibly low computational cost.
\end{itemize}

 In this paper we explore the potential of such an adaptive solution strategy on a primarily
 theoretical level.
 Our central objective is to analyze the intrinsic complexity of accuracy oriented
 eigenvalue  computations.
 In particular, we shall show under which circumstances such a scheme
 has asymptotically optimal complexity in the following sense:
 If the solution $u$ belongs to some approximation space $\cA^s$, which means
 that $N$ terms in the expansion of $u$ suffice to approximate $u$ in $\cH$ within
 accuracy of order $N^{-s}$, then the computational work for realizing a target $\ve$
 remains bounded by a fixed multiple of $\ve^{-1/s}$, when $\ve$ tends to zero.

 While the present work is certainly inspired by prior work on adaptive solution schemes for
 operator equations (see e.g. \cite{CDD2,CDD3}), there are some noteworthy differences.
 On one hand the iteration is nonlinear, on the other hand, one approximates coefficient
 arrays that are only determined up to normalization.

 The layout of the paper is as follows. After summarizing the fundamental properties of the problem in Section \ref{sectproblem}, we formulate in Section \ref{sect2} an iterative scheme
 for the computation of the smallest eigenvalue of the infinite dimensional problem
  \eref{1.0} (or equivalently \eref{1.4a})
 and a corresponding eigenspace. Moreover, we analyze the convergence of this idealized iteration
 making essential use of the ellipticity of $\bA-\lambda\bC$ on certain orthogonal complements. Section \ref{adapsec} is devoted to the numerical realization
 of the idealized scheme based on the approximate realization of residuals. We first show under which
 circumstances and for which accuracy tolerances such perturbed schemes converge to
 the exact solution of \eref{1.4a}. Then, in a second step, we analyze the complexity of such schemes
 under the assumption that the involved operators $\bA, \bC$ are in a certain sense {\em quasi-sparse}.
 We note that this property is known to hold for a wide range of operators. The main result is
 that the proposed adaptive scheme exhibits in some sense asymptotically optimal complexity.
 While the underlying adaptive scheme is to be viewed as one possible prototype realization
 with certain asymptotic properties, we conclude the section with indicating ways of
 quantitative improvements.
 Finally, in Section \ref{schroedinger} we briefly indicate a typical application background.

\section{Problem formulation}\label{sectproblem}

Motivated by the above considerations, we will in the remainder of this paper restrict our treatment to determining the smallest eigenvalue $\underline \lambda$ and the associated eigenvector $\mathbf{u}$
of the generalized eigenvalue problem
\beqn
 \label{evp}
\mathbf{A}\bu= \underline \lambda {\mathbf C}\bu,
\eeqn
formulated on the infinite dimensional
space $\ell_2(\cI)$. As in the previous examples we will only consider the case where $\mathbf{A, C}$ are symmetric, positive definite matrices defined on $\ell_2(\cI)$ endowed as above with the norm $\|\bv\|^2:=\sum_{i\in \cI}
|x_i|^2=: \lll \bv,\bv\rr$.
Let us again stress that the index set $\cI$ could (and actually will) be infinite.

Further, we will impose the following conditions on our problem \eref{evp}
reflecting the properties we have identified in the framework discussed before.

\begin{assumption}
There exist positive constants $\gamma, \Gamma$ such that
\beqn
~ \qquad \gamma ~||{\bx}||^2 ~~\leq~ ~||{\bx}||^2_{\mathbf A}~ ~\leq~ ~\Gamma ~||{\bx}||^2,\label{Aund2}
\eeqn
where $||\cdot ||_{\mathbf A} := \langle {\mathbf A}\cdot ,.\cdot \rangle^{\frac{1}{2}}$ denotes the ${\mathbf A}$-(energy)-norm, see \eref{1.4c}.
\end{assumption}

\begin{assumption}
\label{raylass} For the minimal generalized simple eigenvalue
$\underline \lambda$ of \ref{evp}, there holds \beqn
\label{raylminl2} 0 ~~< ~~\underline \lambda ~=~ \min_{\bx}
\frac{\lll \bA\bx,\bx\rr}{\lll \bC\bx,\bx\rr}, \eeqn while there
exists a $\Lambda > \underline{\lambda}$ one has for all $\bx$ with
$\langle \bA \bu, \bx \rangle = 0$ \beqnr \label{orthgapl2}
\frac{\langle \bA \bx, \bx \rangle}{\lll \bC\bx,\bx\rr} ~\geq
~\Lambda ~>~ \underline \lambda. \eeqnr
 \end{assumption}

\begin{remark}
In the framework discussed in the previous section it was already observed
that $\bC$ is bounded. Note that this
actually follows from  Properties 1, 2, namely
one has for  $C_\bC:= \Gamma/\underline \lambda$ that
\beqn
~ \qquad ~||{\bx}||_{\mathbf C}~ ~\leq~ ~C_{\mathbf C} ~||{\bx}||.
\label{Cund2}
\eeqn
On the other hand, it follows from \eref{Aund2} together with \eref{raylmin} that $\bC$ is coercive on all of $\ell_2$ if and only if $\cL$ is bounded with respect to the $\XX$-inner product, a condition that is usually not fulfilled, for instance
when $\cL$ is a differential operator.
\end{remark}
%

\section{Basic Algorithm - A Richardson-style method for calculating the smallest eigenvalue}
\label{sect2}

In this section, we shall now formulate an ``ideal iteration
scheme'' for \eref{evp}. This scheme and the tools used to prove
convergence will be utilized later in the analysis for our adaptive
algorithm to be introduced in Section \ref{adapsec}; therefore, the
convergence analysis will be given in detail. Let $\underline{\tilde
\lambda}$ and $\tilde \Lambda$ be lower and upper bounds for
$\underline\lambda$ and $\Lambda$, respectively. Then the iteration
reads as follows.
\subsection{Basic Algorithm:}\label{Richie}

\bigskip
\textit{\textbf{MINIT}}\vspace*{0.5mm} 
\footnotesize
\hrule \medskip\noindent
\textit{\textbf{Require:} initial value} ~$\mathbf{x}_0 \in \ell_2 (\mathcal{I}), \|\mathbf{x}_0\|=1, \\ \alpha = 2 ((1 - \underline{\tilde \lambda}/\tilde \Lambda ) \gamma + \Gamma)^{-1}$\\
\textit{\textbf{Iteration:} \\
For $n= 0, 1, \ldots$ do\\
\phantom{ooo}Calculate the Rayleigh quotient $\lambda^{(n)}= \frac{\langle \mathbf{Ax}_n,\mathbf{x}_n \rangle}{\langle \mathbf{Cx}_n,\mathbf{x}_n\rangle};$\\
\phantom{ooo}Calculate the residual $\mathbf{r}_{n}= \mathbf{A x}_n-\lambda^{(n)} \mathbf{C x}_n.$\\
\phantom{ooo}Let $\widehat{\mathbf{x}}_{n+1} = \mathbf{x}_n - \alpha \mathbf{r}_{n} $;\\
\phantom{ooo}Normalize $\mathbf{x}_{n+1}=\langle \widehat{\mathbf{x}}_{n+1},\widehat{\mathbf{x}}_{n+1}\rangle^{-\frac{1}{2}}\widehat{\mathbf{x}}_{n+1}.$\\
endfor} \medskip
\hrule
\normalsize
\medskip
This algorithm resembles a preconditioned Richardson  iteration
$\widehat{\mathbf{x}}_{n+1} = \mathbf{\Phi}^{(n)} \mathbf{x}_n$ with
a stage dependent iteration matrix
$$\mathbf{\Phi}^{(n)}:=
\mathbf{I}- \alpha ( \mathbf{A} - \lambda^{(n)} \mathbf{C})$$
depending on the $n$-th iterate $\mathbf{x}_n$ and a relaxation
parameter $\alpha$. On the other hand, defining the (generalized)
 Rayleigh quotient $\lambda(\mathbf{x})=\frac{\langle \mathbf{Ax},\mathbf{x}\rangle}{\langle \mathbf{Cx},\mathbf{x}\rangle}$  one has  $\nabla \lambda(\mathbf{x}_n) = \frac{2}{\langle \mathbf{Cx}_n,\mathbf{x}_n \rangle}\mathbf{r}_{n}$, so that a steepest descent step with stepsize $\beta_n$ takes the form
 $$
 \bx_n -\frac{2\beta_n}{\lll \bC\bx_n,\bx_n\rr} \mathbf{r}_n.
 $$
 Thus, the above iteration could be viewed as a steepest descent method for the generalized
 Rayleigh quotient. This iteration scheme is also known as a
 preconditioned inverse iteration scheme
 PINVIT, see e.g. \cite{knyazev-ney,rsz}. Note that the preconditioning is already provided through the formulation as a well-conditioned problem in $\ell_2(\mathcal{I})$.

This algorithm can be improved by subspace acceleration techniques.
For instance, the Rayleigh quotient might be minimized according to \beqn
\label{linesearch} \bx_{n+1}:= {\rm argmin}\,\{\lambda(\bvv) :
\bvv\in {\rm span}\,\{\bx_n,\br_n\}, ||\bvv||=1\}. \eeqn
This corresponds to optimal line search in a steepest descent step. However,
in order to have a technically simpler exposition of the principal mechanisms
of our approach we confine the subsequent discussion first to the
above simple Richardson type iteration and will take up
possible improvements along with computational consequences later in Section \ref{sect3.3}.

In the remainder of this section, we will analyze the properties of
Algorithm \textit{\textbf{MINIT}} and in particular prove
convergence to the smallest eigenvector, provided the starting
vector lies sufficiently close to the target $\mathbf u$. The main
results will then be compiled in Theorem \ref{Richkonvergenz} at the
end of this section.

\subsection{Error reduction} \label{richkonvsec}
We shall show that for suitable damping parameters $\alpha$ the above (ideal) iteration
(on the full infinite dimensional space $\ell_2(\cI)$)
provides approximations to the searched eigenpair with a guaranteed fixed error reduction per step.
\subsubsection{Preliminary considerations}

Setting $\mathbf{\boldsymbol{\delta}}_n
 :=\mathbf{\boldsymbol{\delta}}(\mathbf{x}_{n}):=\mathbf{x}_n - \mathbf{u}$
for the normalized iterates $\mathbf{x}_n$, we will  show that the
error components that are orthogonal to $\mathbf{u}$ tend to zero in
the Euclidean $\ell_2(\cI)$-norm. Denoting by $\mathbf{P}\bx =
\frac{\langle \mathbf{x}, \mathbf{u} \rangle }{\langle \mathbf{u},
\mathbf{u} \rangle}\mathbf{u}$  the $\ell_2$-orthogonal projection
of $\bx$ into the eigenspace $U_0$ containing the eigenvector $\bu$,
the orthogonal error component is given by $
\boldsymbol{\delta}_{n}^{\bot} = (\mathbf{I} - \mathbf{P})
\mathbf{\boldsymbol{\delta}}_{n} $, we introduce
 \beqn
\|\boldsymbol{\delta}_{n}^{\bot}\| = || (\mathbf{I} - \mathbf{P})
\mathbf{\boldsymbol{\delta}}_{n}||= \| (\mathbf{I} -
\mathbf{P} )\mathbf{x}_n \| ~=:~ \sin \angle(\mathbf{x}_n,
\mathbf{u})
 \ .
 \eeqn
 Note that $\|\boldsymbol{\delta}_{n}^{\bot}\|$ does not depend on the normalization of $\mathbf{u}$.
We shall now show that we have an at least linear reduction of the orthogonal error component
of the iterates $\mathbf{x_{n}}$ of \eref{Richie}, i.e.
\beqn
 \label{linkon}
  ||\mathbf{\boldsymbol{\delta}}^{\bot}_{n+1}||\leq \xi \, ||\mathbf{\boldsymbol{\delta}}^{\bot}_{n}||,
   \eeqn
   with $\xi < 1,$ provided the initial vector has a sufficiently small angle with
    the solution $\mathbf{u}$ corresponding to the smallest eigenvalue $\underline \lambda$.

Towards this end, we first note that $\mathbf{x}_n \bot~\mathbf{r}_n$, so that we have for the non-normalized iterates
\beqn
\label{2.8}
|| \widehat{\mathbf{x}}_{n+1}  ||^2= || \mathbf{x}_{n} - \alpha \mathbf{r}_n ||^2 = || \mathbf{x}_{n}||^2  + \alpha^2 || \mathbf{r}_n ||^2 \geq ||\mathbf{x}_{n}||^2 = 1.
 \eeqn
Hence,  the error of the $\ell_2$-normalized iterates $\mathbf{x}_{n+1}$ can be estimated by
\beqn
||\mathbf{\boldsymbol{\delta}}^{\bot}_{n+1}||~  =~ \frac{1}{||\widehat{\mathbf{x}}_{n+1}||} ||(\mathbf{I} - \mathbf{P}) \widehat{\mathbf{x}}_{n+1}||~ \leq~ ||(\mathbf{I} - \mathbf{P}) \widehat{\mathbf{x}}_{n+1}|| ~=:~ ||\widehat{\mathbf{\boldsymbol{\delta}}}_{n+1}^{\bot}||.
\label{normmachtklein}
\eeqn
To go further, we decompose  $\widehat{\mathbf{\boldsymbol{\delta}}}_{n+1}^{\bot}$ as follows
 \begin{eqnarray}
 \label{2.10}
  ||\widehat{\mathbf{\boldsymbol{\delta}}}_{n+1}^{\bot}||&=&
   ||(\mathbf{I} - \mathbf{P})\mathbf{\Phi}^{(n)}\mathbf{x}_n||  
= ||(\mathbf{I} - \mathbf{P})\left(~ (\mathbf{I} - \alpha ( \mathbf{A} - \lambda^{(n)}\mathbf{C})\right)\mathbf{x}_n || \nonumber\\
&\leq& ||(\mathbf{I} - \mathbf{P})~(\mathbf{I} - \alpha ( \mathbf{A} - \underline \lambda \mathbf{C} )) \mathbf{\boldsymbol{\delta}}_n ~ ||+  ||\alpha(\lambda^{(n)} - \underline \lambda) (\mathbf{I} - \mathbf{P})\mathbf{ C x}_n||,
\end{eqnarray}
where we have used that
$$\bM\bu :=(\bA -\underline\lambda \bC)\bu=\mathbf{0}\ .$$
 Moreover, using the fact that
${\rm range} {\mathbf{P}} \subseteq \ker \mathbf{M}$ so that
 \beqn
  \mathbf{M} ~=~ \mathbf{M} (\mathbf{I}-\mathbf{P}),
  \label{projektoreinfuegen}
  \eeqn
we conclude, upon using $ (\mathbf{I}-\mathbf{P})^2= \mathbf{I}-\mathbf{P}$ and
$$
( \mathbf{A} - \underline\lambda \mathbf{C} )) \mathbf{\boldsymbol{\delta}}_n=\bM  \mathbf{\boldsymbol{\delta}}_n=
\bM (\mathbf{I}-\mathbf{P})  \mathbf{\boldsymbol{\delta}}_n =\bM \mathbf{\boldsymbol{\delta}}_n^\perp,
$$
  that
\beqn
||\widehat{\mathbf{\boldsymbol{\delta}}}_{n+1}^{\bot}|| ~~\leq~~
  ||~(\mathbf{I} - \mathbf{P})~\left(\mathbf{I} - \alpha \mathbf{M} \right)
   \mathbf{\boldsymbol{\delta}}^{\bot}_n ~ ||+ C_{\bC}~\alpha
   ~|\lambda^{(n)} - \underline \lambda |. ~\label{zerlegt} \eeqn
We now estimate the two parts in the right hand side of
\eref{zerlegt} separately in the next two subsections.

\subsubsection{Linear part of the iteration scheme}

The main goal of this subsection will be to show the following
property of the iteration matrix, from which in Lemma \ref{Rich(a)}, an estimate for the
left part of \eref{zerlegt} will easily be derived.

\begin{lemma} \label{Mell}
Under the assumptions stated above, the matrix  %
 \beqn
 \label{M}
 \mathbf{M} := \mathbf{A} - \underline\lambda \mathbf{C}
 \eeqn
 is bounded and $\ell_2$-elliptic on the set
 \beqn
 \label{V_0}
 V_0 = \{ \bx \in \ell_2 ( \mathcal{I}) : \langle \bx,\mathbf{u} \rangle= 0 \},
 \eeqn
that is, there exist $\theta, \Theta > 0$ such that
\beqn\label{2.5}
\theta~|| \bx||^2 ~~\leq ~~\langle \mathbf{M} \bx, \bx \rangle ~~\leq~~
\Theta ~||\bx||^2~~~ \makebox{ for all}~~~ \bx \in V_0.
 \eeqn
\end{lemma}

For the proof of Lemma \ref{Mell}, we need the following general facts
included for the convenience of the reader in form of the next two lemmata.

\begin{lemma}
 \label{projectornorm}
Let $(\XX,\langle\cdot , \cdot \rangle)$ be any inner product space with norm
$||\cdot||^2:= \lll \cdot,\cdot\rr$. Moreover,
assume that  $x, v \in X$ fulfil the conditions $||v|| = ||x||$
and  $\langle x, v\rangle \geq 0$. Then, for  $P = \frac{\langle .,v \rangle}{\langle v,v \rangle} v$,
we have
$$
\frac{1}{\sqrt{2}} ||x -v|| ~\leq ~|| (I-P) x|| ~\leq~ ||x - v ||.
$$
\end{lemma}
\noindent
{\bf Proof:}
First, we note that it suffices to show the equivalence for $||v|| = ||x|| =1$ and $v \neq x$.
The second inequality is clear due to the fact that $P$ is an $\ell_2$-orthogonal projector. For the inverse estimate, let $\alpha_0 = \langle x, v\rangle \geq 0$. We then have $ ||x-v||^2 ~=~  ||x||^2 + ||v||^2 - 2 \alpha_0 ~=~ 2 - 2 \alpha_0 $ and $|| (I-P) x ||^2 = 1- \alpha_0^2,$ which together gives
$$
\frac{ || (I - P) x ||^2 }{||v - x ||^2} ~=~ \frac{ 1 + \alpha_0}{ 2 }~\geq~ \frac{1}{2}.
$$

\feddich
 \begin{lemma}
\label{ellipticity} Let $\langle \cdot, \cdot \rangle_1$, $\langle
\cdot, \cdot \rangle_2$ be two equivalent inner products on a
Hilbert space $\cX$, i.e. they induce equivalent norms in $ \cX$,
 $G$ a symmetric operator with respect to $\langle\cdot ,\cdot \rangle_2$
 (i.e. $\langle Gx, y \rangle_2 = \langle x, Gy \rangle_2$ for all $x, y \in \cX$) and $u \in \ker G$. Then, if $G$ is $\lll\cdot,\cdot\rr_2$-elliptic on the $\lll\cdot,\cdot\rr_1$-orthogonal complement of $u$, that is
\beqn \label{compell} \lll Gx,x\rr_2\geq c \lll x,x\rr_2\eeqn holds for all $x$ with $\lll x,u\rr_{1} = 0,$ then
$G$ is also $\lll\cdot,\cdot\rr_2$-elliptic on the $\lll\cdot,\cdot\rr_2$-orthogonal complement of $u$, i.e. \eref{compell} holds for all $x$ with $\lll x,u\rr_{2} = 0$ with a possibly different constant $c$.
\end{lemma}
\noindent
{\bf Proof:} ~
Let $x \in \cX$ fulfil $\langle x, u \rangle_2 = 0$. Decomposing $x = P_1 x + x^{\bot_1}$, where $P_i x= \frac{\langle x,u\rangle_i}{ \langle u,u\rangle_i} u$ denotes the $i$-orthogonal projector, we obtain
$Gx = Gx^{\bot_1}$ and therefore
$$
\langle Gx, x \rangle_2 ~=~ \langle x^{\bot_1}, Gx \rangle_2 ~=~\langle Gx^{\bot_1}, x^{\bot_1} \rangle_2 ~\gsim  || x^{\bot_1}||_2^2~\sim~||x^{\bot_1}||_1^2.
$$
To prove the assertion, it remains to see that $|| x^{\bot_1}||_1
\sim || x^{\bot_2}||_2 = || x ||_2$, where $x^{\bot_2} = x -
\frac{\langle x, u\rangle_2 }{\langle u, u\rangle_2 }u = x$. This
latter fact follows from Lemma \ref{projectornorm} above, where we
can choose $v$ to be a multiple of $u$ such that $||v||_1=||x||_1$
and $\langle x, v \rangle_1 \geq 0$, giving
$$
 || x^{\bot_1}||_1 ~ \sim ~|| x - v||_1 ~\sim || x - v||_2 \geq ||(I-P_2)( x - v)||_2 = || x^{\bot_2}||_2.
 $$
 This confirms the assertion.
\feddich

\noindent
{\bf Proof of lemma \ref{Mell}.} For the ellipticity, we show the claim for all $\bx$ with $\lll \bx, \bu \rr_{\bA} = 0$; setting $\lll \cdot,\cdot\rr_1 = \lll \cdot,\cdot\rr_{\bA}$ in Lemma \ref{ellipticity} proves then the coercivity on $V_0$. Indeed, for all such $\bx$, we have $\frac{||\bx||_{\bA}^2}{||\bx||_{\bC}^2} \geq \Lambda$, by assumption \ref{raylass}, and therefore
\begin{align} \lll (\bA-\underline\lambda\bC)\bx,\bx\rr ~&= ~\| \bx\|_{\bA}^2 - \underline \lambda \|\bx\|_{\bC}^2 \\
   &= ~\frac{\Lambda - \underline \lambda}{\Lambda}\| \bx\|_{\bA}^2  + \frac{\underline \lambda}{\Lambda} \| \bx\|_{\bA}^2 - \underline \lambda \|\bx\|_{\bC}^2 \\
   &\ge ~\frac{\Lambda - \underline \lambda}{\Lambda} \gamma \| \bx\|^2 +   \frac{\underline \lambda}{\Lambda} \Lambda \| \bx\|_{\bC}^2 - \underline \lambda \|\bx\|_{\bC}^2  \\
   & =~ \frac{\Lambda - \underline \lambda}{\Lambda} \gamma \| \bx\|^2.
\end{align}
 By
\begin{align}
 0 ~<  ~ \lll (\bA-\lambda\bC)\bw,\bw\rr~ \le ~\lll  \bA\bw,\bw\rr \le \Gamma \|\bw\|^2,
\end{align}
the boundedness of $\bM$ on $V_0$ follows.

\feddich

\begin{lemma}
\label{Rich(a)}
Let $\mathbf{\Phi} :=(\mathbf{I} - \mathbf{P})(\mathbf{I} - \alpha  \mathbf{M}): V_0 \rightarrow V_0$
denote the matrix in the left part of (\ref{zerlegt}).
Then setting the parameter $\alpha := 2 ((1 - \underline \lambda/ \Lambda) \gamma + \Gamma)^{-1}$
in the Richardson method gives
$$
||\mathbf{\Phi }||_{V_0 \to V_0}~ \leq~ \beta~ < ~1.
$$

\end{lemma}
\noindent
{\bf Proof:}
Due to the continuity and ellipticity of $\mathbf{M}$ on $V_0$, we have
\beqn
 (1 - \alpha \Theta ) ||\mathbf{x}||^2 ~\leq~ \langle(\mathbf{I} - \alpha  \mathbf{M}) ~\mathbf{x}, \mathbf{x} \rangle ~ \leq~ (1 - \alpha \theta ) ||\mathbf{x}||^2,~~~\forall~~\mathbf{x} \in V_0.
 \label{ellip1}
 \eeqn
   The same inequality applies to $\mathbf{\Phi}$ owing to the fact that $\mathbf{P}$
   is symmetric and therefore
$$ \langle (\mathbf{I} - \mathbf{P})(\mathbf{I} - \alpha  \mathbf{M})
 \mathbf{x}, \mathbf{x}\rangle = \langle (\mathbf{I} - \alpha  \mathbf{M})
  \mathbf{x}, (\mathbf{I} - \mathbf{P}) \mathbf{x}\rangle = \langle (\mathbf{I} - \alpha
   \mathbf{M}) \mathbf{x}, \mathbf{x}\rangle,~~~\forall~~\mathbf{x} \in V_0.
$$
 We now let $ \beta = \frac {\Theta - \theta}{\Theta + \theta} < 1$. Then,
 choosing $\alpha := \frac{2}{\Theta + \theta} = 2 ((1 - \underline \lambda/
  \Lambda) \gamma + \Gamma)^{-1}$ (recalling the choice of $\Theta$
   and $\theta$ from the proof of Lemma \ref{Mell}), equation (\ref{ellip1})  gives
$-\beta \leq \langle \mathbf{\Phi} ~\mathbf{x} ,\mathbf{x}  \rangle
\leq \beta$ and therefore $ ||\mathbf{\Phi}||_{V_0 \to V_0} \leq
\beta < 1$. \feddich

\subsubsection{Estimates for Rayleigh quotients and residuals}
It remains to estimate the Rayleigh quotients in the right hand side of equation (\ref{zerlegt}).

\begin{lemma}
\label{Rich(b)} Let $\mathbf{x} \in \ell_2(\cI)$ with $||\mathbf{x}
|| = 1$ and as before let  $\lambda(\mathbf{x})$ denote the
corresponding Rayleigh quotient. Then we have \beqnr
\lambda(\mathbf{x})- \underline\lambda~  \quad &\leq&  \quad
\frac{2\Gamma}{\gamma} ~\lambda(\mathbf{x})~
||\mathbf{\boldsymbol{\delta}}(\bx)^{\bot}||^2,
 \label{ewzeile2}
 \eeqnr
where the constants in inequality (\ref{ewzeile2}) are those from the norm estimate (\ref{Aund2}) and
$\mathbf{\boldsymbol{\delta}}^{\bot}(\bx)=(\mathbf{I}-\mathbf{P})(\bx -\bu)$ as above.
\end{lemma}
\noindent
{\bf Proof:}
Suppose that
 $\mathbf{u}_x$ is an  eigenvector belonging to the smallest eigenvalue $\underline \lambda$  of the eigenvalue problem (\ref{evp}) which is normalized in such a way that
  $||\mathbf{u}_x|| = ||\mathbf{x}||$
  and $\langle \mathbf{x}, \mathbf{u}_x \rangle \geq 0$.
  Straightforward computation yields
  \beqnr \label{ewAC}
\lambda(\mathbf{x})- \underline\lambda~ \quad =
 \quad \frac{1}{||\mathbf{x}||_\mathbf{C}^2}\left(
 ||\mathbf{x}-\mathbf{u}_x||_ {\mathbf{A}}^2-
 \underline\lambda ||\mathbf{x}- \mathbf{u}_x||_{\mathbf{C}}^2\right). \eeqnr

Keeping the normalization of $\bx$ and \eref{Aund2} in mind and
since there are no absolute values  involved there, we infer from \eref{ewAC} that
\beq 
\lambda(\mathbf{x})- \underline \lambda & \leq & \frac{1}{||\mathbf{x}||^2_{\mathbf{C}}} ||\mathbf{x}-\mathbf{u}_x||_ {\mathbf{A}}^2 ~=~  \frac{||\mathbf{x}||^2}{||\mathbf{x}||^2_{\mathbf{C}}} ||\mathbf{x}-\mathbf{u}_x||_ {\mathbf{A}}^2~ \leq ~ \frac{||\mathbf{x}||_\mathbf{A}^2}{\gamma ||\mathbf{x}||^2_{\mathbf{C}}} ||\mathbf{x}-\mathbf{u}_x||_ {\mathbf{A}}^2~ \\ &=& ~\frac{\lambda(\mathbf{x})}{\gamma}  ||\mathbf{x}-\mathbf{u}_x||_ {\mathbf{A}}^2 ~ ~~\leq ~~ \frac{\Gamma}{\gamma} \lambda(\mathbf{x}) ||\mathbf{x}-\mathbf{u}_x||^2 .
\eeq
Since Lemma \ref{projectornorm} says that
$$
|| \bx -\bu_x|| \leq \sqrt{2}||(\mathbf{I} -\mathbf{P})\bx|| =\sqrt{2}||(\mathbf{I} -\mathbf{P})(\bx-\bu_x)||,
$$
the assertion \eref{ewzeile2} follows as well.
\feddich

\begin{remark}
\label{rem2.1}
As before let $\boldsymbol{\delta}=\boldsymbol{\delta}(\bx):=\bx-\bu$ and fix $0<a <\gamma/2\Gamma$. If
$||\boldsymbol{\delta}(\bx)^\bot||\leq \sqrt{a}$, one has for the corresponding Rayleigh quotient
\beqn
\label{raylabs}
\lambda(\bx)\leq \underline \lambda\big(1-\frac{2\Gamma a}{\gamma}\big)^{-1}=: K.
\eeqn
Taking e.g. $a=\gamma/4\Gamma$ one has $K=2\underline \lambda$.
\end{remark}
\noindent
{\bf Proof:} From Lemma \ref{Rich(b)} we infer that
$\lambda(\bx)\leq \underline \lambda\big(1-\frac{2\Gamma }{\gamma}||\boldsymbol{\delta}(\bx)^\bot||^2\big)^{-1}
\leq \underline \lambda\big(1-\frac{2\Gamma a}{\gamma}\big)^{-1}$.
\feddich

An immediate consequence of  \eref{raylabs} can be stated as follows.
\begin{remark}
\label{remCcoercive}
Under the property of Section \ref{sectproblem}  the operator $\bC$ is coercive on any sufficiently small fixed angular neighborhood of
the direction $\bu$, i.e. for $K$ from \eref{raylabs} one has 
\beqn
\label{Ccoercive}
\lll\bC\bx,\bx\rr \geq \gamma/K \quad
\mbox{whenever} \,\,||\bx|| =1, \,\, ||(\bu-\bx)^\bot ||\leq \sqrt{a}.
\eeqn
\end{remark}

The next observation is that the orthogonal error components are controlled by the residuals.
To this end,  note that
$$
\mathbf{M}_{\bot}:~~ V_0 \rightarrow V_0,~~ \mathbf{x} ~~\mapsto~~ (\mathbf{I}-\mathbf{P})\mathbf{Mx}
$$
 is bounded and elliptic so that $||M_\bot^{-1}||$ is a bounded positive number.

\begin{lemma}
\label{lem:lema}
For any $\bx\in \ell_2(\cI)$ with $||\bx||=1$, let
$  \mathbf{r}(\bx) :=   ( \mathbf{A} - \lambda(\bx) \mathbf{C})\mathbf{x}$ be the corresponding residual with respect to the Rayleigh quotient $\lambda(\bx):= \frac{\langle  \mathbf{A} \mathbf{x} ,
\mathbf{x}\rangle}{\langle \mathbf{C}\mathbf{x},
 \mathbf{x}\rangle } $ and, as before, let   $\mathbf{x}^{\bot}=(\mathbf{I- P})\bx$  so that in particular
 $\boldsymbol{\delta}(\bx)^{\bot}=(\mathbf{I-P})(\bx -\bu)= (\mathbf{I-P})\bx$. Assume that for the constant
 $a$ from Remark \ref{rem2.1}
\beqn
\label{deltaconda}
||\boldsymbol{\delta}(\bx)^\bot || \leq \min\,\Big\{ \sqrt{a}, \frac 12\Big\{ \Big(1+\frac{\gamma}{\Gamma K C_{\bC}
 ||\bM_\bot^{-1}||}\Big)^{1/2}-1\Big\}\Big\}.
\eeqn
Then for $M:= 2||\bM_\bot^{-1}|| $ one has
\beqn
\label{deltares}
\|\boldsymbol{\delta}(\bx)^{\bot} \| ~~\leq~~  M \|\mathbf{r}(\bx)^{\bot} \| ~~\leq~~ M||\mathbf{r}(\bx)|| ~~\leq ~~
M\Big(||\bM|| + \frac{2K\Gamma C_{\bC}}{\gamma}||\boldsymbol{\delta}(\bx)^\bot||\Big) ||\boldsymbol{\delta}(\bx)^\bot||.
\eeqn
\end{lemma}

\noindent
{\bf Proof:}
Since
 $$
 \mathbf{r}(\bx) = \bM(\bx -\bu)+(\underline \lambda -\lambda(\bx))\bC\bx=
 \bM(\mathbf{I-P})(\bx-\bu)+(\underline \lambda -\lambda(\bx))\bC\bx,
 $$
 we conclude, on account of Lemma \ref{Rich(b)} and Remark \ref{rem2.1},  that
\beqn
||\mathbf{r}(\bx)|| \leq ||\bM||\, ||\boldsymbol{\delta}^\bot(\bx)|| + \frac{2K\Gamma C_{\bC}}{\gamma}||\boldsymbol{\delta}^\bot(\bx)||^2,
\eeqn
which confirms the upper estimate in \eref{deltares}.

As for the lower estimate, recall that by the definition of $\bM_\bot$,
$$
\bM_\bot \boldsymbol{\delta}(\bx)^\bot =  (\mathbf{I}-\mathbf{P}) (\mathbf{A}- \lambda(\bx) \mathbf{C}  )
 \boldsymbol{\delta}(\bx)^\bot
+  (\mathbf{I}-\mathbf{P}) (\lambda(\bx)- \underline \lambda) \mathbf{C}
\boldsymbol{\delta}(\bx)^{\bot},
$$
 giving
 \begin{equation}
 \label{2.23}
\|\mathbf{M}_{\perp }^{-1} \|^{-1}
\boldsymbol{\delta}(\bx)^{\bot}\|   \leq   \|\mathbf{M}_{\perp } 
\boldsymbol{\delta}(\bx)^{\bot} \|
 \leq
\| (\mathbf{I}-\mathbf{P}) (\mathbf{A}- \lambda(\bx) \mathbf{C}  )  
\boldsymbol{\delta}(\bx)^{\bot} \| ~+ ~\|   (\lambda ( \mathbf{x}) -
\underline \lambda )
\mathbf{C}
\boldsymbol{\delta}(\bx)^{\bot} \|.
\end{equation}
Now straightforward calculations yield
$$
 (\mathbf{A}- \lambda(\bx) \mathbf{C}  ) \boldsymbol{\delta}(\bx)^{\bot}=\mathbf{r}(\bx)- (\mathbf{A}- \lambda(\bx) \mathbf{C}  )
 \mathbf{P}\bx = \mathbf{r}(\bx)- \frac{\lll \bx,\bu\rr}{\lll \bu,\bu\rr}(\underline \lambda -\lambda(\bx))\bC\bu,
 $$
which gives
\beqn
\label{2.24}
\| (\mathbf{I}-\mathbf{P}) (\mathbf{A}- \lambda(\bx) \mathbf{C}  )
\boldsymbol{\delta}(\bx)^{\bot} \| \leq ||\mathbf{r}(\bx)^\bot || + C_{\bC} \,|\underline \lambda -\lambda(\bx)|.
\eeqn
Invoking now Lemma \ref{Rich(b)}, \eref{ewzeile2} together with Remark \ref{rem2.1} and \eref{2.23},
yields
\beqn
\label{2.25}
\|\mathbf{M}_{\perp }^{-1} \|^{-1}
\|\boldsymbol{\delta}(\bx)^{\bot}\|   \leq ||\mathbf{r}(\bx)^\bot || + \frac{2\Gamma K C_{\bC} }{\gamma}\big( 1+
||\boldsymbol{\delta}(\bx)^{\bot}||\big)
||\boldsymbol{\delta}(\bx)^\bot||^2.
\eeqn
Hence, whenever
\beqn
\label{2.26}
 \frac{2\Gamma K C_{\bC}}{\gamma}\big( 1+||\boldsymbol{\delta}(\bx)^{\bot}||\big)
||\boldsymbol{\delta}(\bx)^\bot||\leq \frac 12 \|\mathbf{M}_{\perp }^{-1} \|^{-1},
\eeqn
which is the case when \eref{deltaconda} holds,
 we obtain
\beqn
\label{2.27}
\|\boldsymbol{\delta}(\bx)^{\bot}\| \leq 2 \|\mathbf{M}_{\perp }^{-1} \|\, ||\mathbf{r}(\bx)^\bot||,
\eeqn
which is lower estimate in \eref{deltares} with $M:= 2 \|\mathbf{M}_{\perp }^{-1} \|$.
 \feddich

\subsubsection{Convergence of the scheme}

In summary, the above observations allow us to establish the following local convergence
properties of the Richardson eigenvalue iteration.
\begin{theorem} 
\label{Richkonvergenz}
Let the damping parameter $\alpha$ in the basic algorithm \textit{\textbf{MINIT}} from Section
 \ref{Richie} be chosen according  to Lemma \ref{Rich(a)}, so that the bound $\beta$ for
 $||\Phi||_{V_0\to V_0}$ in Lemma \ref{Rich(a)} satisfies $\beta <1$. Furthermore assume
 that the initial error satisfies
 \beqn
\label{deltacond}
||\boldsymbol{\delta}_0^{\bot}||\leq \min\,\Big\{a^{1/2}, \frac{1-\beta}{2\tilde{c} K} 
 \Big\} =: d,
\eeqn
where $\tilde c:=2 \alpha\Gamma C_{\bC}/\gamma$.
 Then the following statements hold:\\[2mm]
 a) For $\xi := (\beta +1)/2<1$ one has 
\beqn
\label{mondelta}
 ||\mathbf{\boldsymbol{\delta}}^{\bot}_{n+1}||\leq \xi \, ||\mathbf{\boldsymbol{\delta}}^{\bot}_{n}||,
  \eeqn
i.e. one has monotone linear error reduction with respect to the $\ell_2$-norm.
Moreover, for any $\epsilon >0$ there is an $n_\epsilon$ such that for
$n\geq n_\epsilon$ one can choose $\xi<\beta +\epsilon$ in \eref{mondelta}, i.e.
one has an asymptotic error reduction of a rate $\leq \beta$.
The orthogonal error components $\mathbf{\boldsymbol{\delta}}_n^{\bot}$ also converge to zero
in the $\mathbf{A}$-norm. \\[2mm]
b)
There exists a uniform constant $C$ such that  the error of the Rayleigh quotients $\lambda(\mathbf{x}_n)$  is bounded by
\beqnr
 \lambda(\mathbf{x}_n)-  \underline \lambda   ~\leq~ C ||\mathbf{\boldsymbol{\delta}}_n^{\bot}||^2 ~ \rightarrow ~ 0.\label{evkonv}
\eeqnr
\vspace*{2mm}
c)
There exist  constants
 $ M, M'>0 $ and $n_0\in \N$ such that
\beqn
\label{deltaresn}
\|\boldsymbol{\delta}_n^{\bot} \|~~ \leq ~~ M \|\mathbf{r}_n^{\bot} \| ~~\leq~~ M||\mathbf{r}_n||~~ \leq~~ M' ||\boldsymbol{\delta}_n^\bot||,
~~~n\geq n_0.
\eeqn
  Moreover, there exists $ \zeta <1$ and $n_1\in \N$ such that
  \beqn
  \label{resred}
  || \mathbf{r}_{n+1}|| ~~\leq~~ \zeta  || \mathbf{r}_{n}||,
   \quad
  n\geq n_1 .  
\eeqn
\end{theorem}

\noindent
{\bf Proof:}~ Recalling  \eref{normmachtklein}, \eref{2.10} 
and employing
Lemmata \ref{Rich(a)} and \ref{Rich(b)} provides  
\begin{eqnarray}
\label{2.18}
 ||\boldsymbol{\delta}_{n+1}^{\bot}||
&\leq &  ||\mathbf{\Phi}||_{V_0 \to V_0}~ ||\boldsymbol{\delta}_n^{\bot}|| +  \alpha \frac{2\Gamma C_{\bC}}{\gamma} \lambda(\mathbf{x}_n)\, ||\boldsymbol{\delta}_{n}^{\bot}||^2\nonumber\\
&\leq&  (\beta +
\tilde{c}~\lambda(\mathbf{x}_n)\,||\boldsymbol{\delta}_{n}^{\bot}||)||\boldsymbol{\delta}_{n}^{\bot}||
=: \xi\,||\boldsymbol{\delta}_{n}^{\bot}||,
\end{eqnarray}
with $\tilde c:=2 \alpha\Gamma C_{\bC}/\gamma$.
It remains to show that $\xi<1$ for $||\boldsymbol{\delta}_{0}^{\bot}||$  sufficiently small. In fact,
we know from Remark \ref{rem2.1} that
 \beqn
  \lambda(\mathbf{x}_0) \leq 
  K,\label{ewabschaetzung}
  \eeqn
provided that $\boldsymbol{\delta}_0=\boldsymbol{\delta}(\bx_0)$ satisfies $||\boldsymbol{\delta}_0^\bot||\leq \sqrt{a}$
(cf. \eref{deltacond}).
Moreover, by \eref{2.18} for $n=0$, we have
$$
||\boldsymbol{\delta}_1^{\bot}||~~\leq~~ (\beta +\tilde{c}K||\boldsymbol{\delta}_0^{\bot}||)||\boldsymbol{\delta}_0^{\bot}||.
$$
Now note that \eref{deltacond}  implies that $\tilde c K||\boldsymbol{\delta}_0^\bot ||
\leq (1-\beta)/2$ so that
\beqn
\label{xibeta}
\beta +\tilde{c}K||\boldsymbol{\delta}_0^{\bot}||~~\leq~~ \frac{1+\beta}2 =: \xi ~<~1,
\eeqn
 which yields
$$
||\boldsymbol{\delta}_1^{\bot}||< ||\boldsymbol{\delta}_0^{\bot}||,~~~ \lambda(\bx_1)\leq K.
$$
One easily shows now with the aid of Lemma \ref{Rich(b)} and \eref{2.18} inductively that
\beqn
\label{monotone}
||\boldsymbol{\delta}_{n+1}^{\bot}||
\leq \xi_n ||\boldsymbol{\delta}_n^{\bot}||,~~ ~~\xi_n:= \beta +\tilde{c}K||\boldsymbol{\delta}_n^{\bot}||\leq \xi<1,~~~~~ \lambda(\bx_n)\leq K,~~~~n\in \N.
\eeqn
This  guarantees monotone convergence of the iteration. It is
obvious from the above arguments that with
$||\boldsymbol{\delta}_{n}^{\bot}||  \rightarrow 0$ the values of
$\xi_n$ drop down to $\beta$, showing the first assertion.

Convergence in the $\mathbf{A}$-norm is now an immediate consequence of
\eref{Aund2} and \eref{mondelta}
provided that $\boldsymbol{\delta}_0$ satisfies \eref{deltacond}.  This proves a).

The convergence of  the Rayleigh quotients $\lambda(\mathbf{x}_n)$ now follows immediately
from  (\ref{ewzeile2}) and \eref{monotone}, giving (\ref{evkonv}), confirming b).

 Concerning c), we have seen above that under the assumption \eref{deltacond}
 the estimates \eref{monotone} hold. Thus, using the rough estimate $a\leq 1$ one
 obtains
\beqn
 \frac{2\Gamma K C_{\bC} }{\gamma}\big( 1+||\boldsymbol{\delta}_n^{\bot}||\big)
||\boldsymbol{\delta}_n^\bot||~~\leq~~ \frac 12 \|\mathbf{M}_{\perp }^{-1} \|^{-1},~~~ n\geq n_0,
\eeqn
where
$$
n_0 := \left\lceil \frac{\log (\gamma/(8\Gamma K C_{\bC}\|\mathbf{M}_{\perp }^{-1} \|)}{\log \xi}\right\rceil.
$$
Hence, for $n\geq n_0$, the hypotheses of Lemma \ref{lem:lema} are satisfied for $\bx=\bx_n$,
which provides the lower estimate
  of \eref{deltaresn} again with
 $M=2||\bM_\bot^{-1}||$ as in Lemma \ref{lem:lema}. The upper estimate is also an immediate
 consequence of Lemma \ref{lem:lema} and \eref{monotone}.

 As for the remaining claim, we  find that
\beq
 ||\mathbf{r}_{n+1}|| &=& ||(\mathbf{A} - \lambda_{n+1} \mathbf{C}) \mathbf{x}_{n+1}|| \\
 &=&
 ||\widehat{\mathbf{x}}_{n+1}||^{-1}~|| (\mathbf{A} - \lambda^{(n+1)} \mathbf{C}) (\mathbf{x}_{n}- \alpha \mathbf{r}_n) ||~  \\
 &=& ||\widehat{\mathbf{x}}_{n+1}||^{-1}~ || (\mathbf{A} - \lambda^{(n)} \mathbf{C}) \mathbf{x}_n \\
 && \quad  -
  \alpha(
  \mathbf{A} - \underline \lambda \mathbf{C}) \mathbf{r}_n+ (\lambda^{(n)} - \lambda^{(n+1)} ) \mathbf{C} \mathbf{x}_n  - \alpha ( \underline \lambda - \lambda^{(n+1)} ) \mathbf{C} \mathbf{r}_n ||.
\eeq
On account of   Lemma  \ref{Rich(b)}, \eref{ewzeile2}, \eref{raylabs}, and the fact that
  $||\widehat{\mathbf{x}}_{n+1}|| \geq 1$ for the unnormalized iterates, we obtain further
\beq
 ||\mathbf{r}_{n+1}||   &\leq& || (\mathbf{I} - \alpha (\mathbf{A} - \underline \lambda \mathbf{C}))|| ~||\mathbf{r}_n||  + \frac{2 C_{\bC}\Gamma K}{\gamma}
   ||\boldsymbol{\delta}_n^{\bot}||^2\big(2 ||\mathbf{x}_n|| + \alpha  ||\mathbf{r}_n||\big)  \\
   &\leq&
\Big(\beta  + \frac{M^2  2C_{\bC}\Gamma K}{\gamma}\big(2 ~+~ \alpha ||\mathbf{r}_n||\big) ||\mathbf{r}_n||\Big) ~||\mathbf{r}_n||,
\eeq
where $M$ is from \eref{deltares} and where we have also used
  Lemma \ref{Rich(a)}, the continuity of $\mathbf{C}$ and the estimate \eref{deltares} just proven before.
 Recall that  we have $\beta < 1$, thus showing the assertion for $||\mathbf{r}_n|| $ small enough.
 The fact that $||\mathbf{r}_n||$ indeed becomes small  follows from the upper estimate
 in \eref{deltares} and \eref{mondelta} as well.
 \feddich

\section{Adaptive strategies for the eigenvalue problem}\label{adapsec}
There are two points that we would like to stress from the start.
First, all the above considerations are independent of the index set
$\cI$ underlying the space $\ell_2(\cI)$ being finite or infinite as
long as our assumptions on $\bA,\bC$ are satisfied. If the basic
Richardson iteration is considered on a fixed finite dimensional
space, i.e. if we consider a fixed discretization of the original
problem \eref{1.0} like \eref{1.0.1}, the proposed scheme may not be
the most favorable one since it offers at best linear error
reduction. Nevertheless, the property \eref{Aund2} reflects the fact that one
has been able to precondition the
problem sufficiently well in the sense that a fixed error reduction
is achieved {\em independently} of the (fixed) size of $\cI$.

 On the other hand, if one is interested in solving the original problem \eref{1.0}
within a desired target accuracy $\ve$ at possibly low cost (which means to understand how the cost
depends on $\ve$ when $\ve$ tends to zero), the game changes completely. Following \cite{CDD2}
one may think of performing the iteration on the infinite eigenvalue problem \eref{evp} which
is still {\em equivalent} to \eref{1.0}. In fact, solving \eref{evp} within target accuracy $\ve$ provides
a solution to \eref{1.0}, due to the norm equivalence \eref{Aund2} and the mapping property
\eref{Lelliptic}, that has up to a constant factor the same accuracy in the (continuous)  energy norm.
Now, of course, the matrices $\bA, \bC$ are infinite so that even when the current iterate $\bx_n$
has finite support, the next iterate $\bx_{n+1}$ cannot be computed exactly. The idea is therefore
to compute each update only approximately within a suitable dynamically varying
accuracy tolerance in such a way that, on one hand, the convergence to the exact solution
is preserved, while on the other hand, the computational  cost of each perturbed iteration is
possibly low.

In the following sections we shall carry out this program based on
the above Richardson type iteration. It will be seen that,
regardless of the order of the iteration, when applied to a fixed
finite dimensional trial space, such an adaptive scheme may perform
at an asymptotically optimal complexity.

We will proceed in two steps. First, each iteration step requires the approximate
calculation of the residual $(\mathbf{A}-\lambda(\bx_n)\bC)\bx_n$ which can be broken down
to two types of tasks, namely:
\begin{itemize}
\item[(I)] to approximate  matrix/vector products $\bA\bx$, $\bC\bx$
(where now $\bA, \bC$ are in principle infinite),
\item[(II)]
to approximate the Rayleigh quotient $\lambda(\bx)= \frac{\lll \bA\bx,\bx\rr}{\lll \bC\bx,\bx\rr}$.
\end{itemize}
Of course, (II) can be reduced to (I) but there are interesting alternatives.
We distinguish the two cases at this point because, depending on the situation, the accuracy requirements may be somewhat different due to
the fact that relative accuracy tolerances are needed in (II) unless $\bC$ has
particularly favorable properties that essentially permit the exact calculation of $\lll \bC\bx,\bx\rr$.
We shall return to this issue later in more detail.

Let us therefore suppose at this point that we have for $\bB\in \{\bA,\bC\}$
a routine with the following property at hand:\\

\noindent
$APPLY(\bB,\bx,\eta)\to \bw$ such that
$$
||\bB\bx -\bw||\leq \eta.
$$

We shall postpone the actual description of $APPLY$ to a later section. Given
such a routine we shall first analyze which tolerances $\eta$ are needed to ensure
convergence of a correspondingly perturbed iteration. As a second step we shall then
discuss the complexity of a corresponding numerical realization.

\subsection {A perturbed iteration} 
\label{adaptalgo}
We shall now collect the main ingredients of a perturbed version of \textit{\textbf{MINIT}} from Section
\ref{Richie}, assuming that the routine $APPLY$ from above is available. The first one concerns the
approximation of a given $ \bx \in \ell_2 ( \mathcal{I})$
 by one with possibly short support. 

\noindent
$APPROX(\bx,\eta)\to \bz$: {\it produces for a given (finitely supported) $\bx$ and any $\eta >0$
a finitely supported $\bz$ such that}
\beqn
\label{3.1.1}
||\bx-\bz||\leq \eta,~~~~\#{\rm supp}\,\bz\leq {\rm argmin}_{||\bw -\bx||\leq \eta/2 }\,\#{\rm supp}\,\bw.
\eeqn

This routine can be realized by replacing as many entries of the input as possible by zero, as long as
the sum of their squares does not exceed $\eta^2$.
One could of course use the tolerance $\eta$ in the $\rm argmin$ which would mean to
compute a best $N$-term approximation of $\bx$. This in turn would require exact sorting of the
coefficients introducing an additional $\log$-factor of the support size of $\bx$, see \cite{CDD1}
for details.
  Being content with quasi-sorting
based on binary binning one can avoid the $\log$-factor at the expense of a slightly larger support,
see e.g. \cite{arne}.

In addition to
approximate matrix/vector products we need approximate Rayleigh quotients
and hence scalar products. Recall that we have to
perform such routines for normalized inputs $||\bx||=1$, which will be henceforth assumed.
A straightforward way to approximate the scalar products would be as follows:\\

\noindent
$SCAL(\bB,\bx,\eta)\to s$:~{\em Given any $\eta>0$, the routine outputs a scalar $s$
  such that
  \beqn
\label{3.1.2}
|\lll \bB\bx,\bx\rr -s|\leq \eta
\eeqn
 as follows:
\begin{itemize}
\item
$APPROX(\bx, \eta/(2||\bB||))\to \bz$;
\item
$APPLY(\bB,\bx,\eta/2)\vert_{{\rm supp}\,\bz} \to \bw$;
\item
$s=\lll \bz,\bw\rr$.
\end{itemize}}

Here $APPLY(\bB,\bx,\eta/2)\vert_{{\rm supp}\,\bz}$ means that the output of
$APPLY$ can be restricted to ${\rm supp}\,\bz$. To which extent this can be used to
economize on the computational effort depends on the realization of $APPLY$ and will
therefore be postponed to Section \ref{sect3.3}.

It is easily seen that $s$ given as above does satisfy \eref{3.1.2}.
 Note that $||\bA ||:=\sup_{||\bx||=1}||\bA\bx|| \leq \Gamma$.

We also remark that one can think of different ways of approximating
$\lll \bB\bx,\bx\rr$. For instance, one could decompose $\bx$ by
determining $\bx_j$ of smallest support so that $||\bx-\bx_j||\leq
2^j\eta^{1/2}$, say, noting that $\bx_J=0$ for $J:= \lceil |\log_2
\eta^{1/2}|\rceil$. Then, setting $\bz_j:= \bx_j-\bx_{j+1}$
\beqn
\label{3.1.3}
\sum_{j=0}^J \lll \bz_j,\bw_j\rr,~~~\bw_j:= APPLY(\bB,\bx,\epsilon_j)\vert_{{\rm supp}\,\bz_j},~~~,j=0\ldots,J,
\eeqn
is an approximation to $\lll \bx,\bB\bx\rr$ of the order $\eta$, provided that $\sum_{j=0}^J
2^j\epsilon_j\sim \eta^{1/2}$. The apparent advantage is that highly accurate matrix applications
need to be computed only on typically small supports of coarse $\bz_j$'s. It will be shown
in Section \ref{sect3.3} that such a strategy does offer asymptotic savings. Of course, as long
as the matrix/vector products $\bA\bx, \bC\bx$ have to be approximated within a similar
accuracy tolerance anyway in the iteration, one might as well stick with the simpler
version of $SCAL$ shown above which we will do for the time being.

Given the routine $SCAL$ we can proceed to approximating the Rayleigh quotients $\lambda(\bx)$.
We shall devise a routine:\\

 \noindent
$RAYL(\bx,\eta)\to \bar\lambda$:  {\it Given $\eta >0$ and any finitely supported $\bx$ with
$||\bx||=1$, $RAYL$ outputs $\bar\lambda$ auch that}
\beqn
\label{3.1.5}
|\lambda(\bx)-\bar\lambda|\leq \eta .
\eeqn

In many cases it is even possible to apply $\bC$ exactly at acceptable cost, e.g. when
$\bC$ is diagonal, see e.g. Remark  \ref{remC}. In this case we need no further assumption on $\bC$ beyond boundedness in $\ell_2(\cI)$
and can simply take
\beqn
\label{3.1.6}
RAYL(\bx,\eta)= \frac{SCAL(\bA,\bx,\eta\lll \bC\bx,\bx\rr)}{\lll\bC\bx,\bx\rr}.
\eeqn

If, on the other hand, $\bC$ can only be applied appoximately with the aid of some routine
$APPLY$, the routine $SCAL$ is slightly more involved. Denoting by $c_C:= \gamma/K$ the local
coercivity constant for $\bC$ from Remark \ref{rem2.1},
 the following realization works.
\begin{lemma}
\label{lemrayl}
Assume that $||(\bu -\bx)^\bot||\leq \sqrt{a}$ holds.
Carrying out the following steps:
\begin{itemize}
\item
$SCAL(\bA,\bx,\eta c_C/2)\to \tilde s_A$,  $SCAL(\bC,\bx, c_C^2\eta/(6\Gamma))\to \tilde s_C$;
\item
$\bar\lambda :=  \tilde s_A/\tilde s_C$,
\end{itemize}
verifies \eref{3.1.5}, provided that
\beqn
\label{3.1.8}
\eta \leq \min\,\{\gamma/c_C, 3\Gamma/c_C\}.
\eeqn
\end{lemma}
{\bf Proof:}
To see
  the validity of \eref{3.1.5}, let $s_A:=\lll\bA\bx,\bx\rr$, $s_C:=\lll\bC\bx,\bx\rr$ and note
first that
$$
\frac{s_A}{s_C}- \frac{\tilde s_A}{\tilde s_C}~= ~\frac{1}{s_C} (s_A-\tilde s_A) -
\frac{\tilde s_A}{s_C\tilde s_C}(s_C-\tilde s_C),
$$
so that
\beqn
\label{3.1.7}
|\lambda(\bx)-\bar\lambda|~\leq ~ \frac{\eta c_C}{2s_C} ~+ ~\frac{\tilde s_A}{s_C\tilde s_C}
\frac{c_C^2\eta}{6\Gamma}.
\eeqn
Now note that under the assumption \eref{3.1.8} one has (recalling that $||\bx||=1$)
$$
|s_A-\tilde s_A|~\leq~ \frac{\eta c_C}{2} ~=~ \frac{\eta c_C||\bx||^2}{2}~\leq~
\frac{\eta c_C s_A}{2\gamma}~\leq~ \frac 12 s_A,
$$
so that $\tilde s_A\leq \frac{3s_A}{2}\leq \frac{3\Gamma}{2}$. Likewise, upon using \eref{Ccoercive},
one obtains for $\eta$ satisfying \eref{3.1.8} that $|s_C -\tilde s_C|\leq s_C/2$ so that
$\tilde s_C\geq s_C/2 \geq c_C/2$. Combining these estimates yields
$$
 \frac{\tilde s_A}{s_C\tilde s_C}~\leq ~\frac{6\Gamma}{2c_C^2},
 $$
which, on account of \eref{3.1.7}, confirms the claim \eref{3.1.5}. \feddich

\begin{remark}
\label{rem3.1.1}
On account of Remark \ref{rem2.1}, the tolerance in the evaluation of the scalar product
remains proportional to the target accuracy $\eta$ in the routine $RAYL$ for either version provided that
$||(\bu -\bx)^\bot||\leq \sqrt{a}$.
\end{remark}

We are now able to describe a routine for approximating the scaled residual
$\alpha(\bA -\lambda(\bx)\bC)\bx$:\\

\noindent
$RES(\bx,\eta)\to \br_\eta$: ~{\it Given any $\eta >0$ and any $\bx$ with $||\bx||=1$, the routine
$RES$ outputs a finitely supported vector $\br_\eta$ such that
\beqn
\label{3.1.9}
||\br_\eta - \alpha(\bA -\lambda(\bx)\bC)\bx||\leq \eta.
\eeqn}
The routine $RES$ can be realized as follows:
\begin{itemize}
\item[1)]
$RAYL(\bx,\eta/(4C_\bC \alpha))\to \bar\lambda$;
\item[2)]
$RES(\bx,\eta)=\alpha\big(APPLY(\bA,\bx,\eta/(2\alpha)) -\bar\lambda APPLY(\bC,\bx,\eta/(4\bar\lambda\alpha))\big)
$
\end{itemize}

Of course, when the approximate evaluation of the scalar product involving $\bA$ simply relies
on the application of $APPLY$ one only needs to carry out this latter routine once in $RES$ with respect
to the minimal tolerance required in $RAYL$ and in step 2) of the above scheme.

\begin{remark}
\label{remparameter} From now on we shall always assume that the
parameter $\alpha$ in the basic Richardson iteration is chosen
according to Lemma \ref{Rich(a)}, ensuring that $\beta< 1$.
\end{remark}

The following lemma will help to understand the perturbed Richarson iteration.
\begin{lemma}
\label{lem3.1.2} Assume that the approximation $\bar\bu$ to $\bu$
satisfies $||\boldsymbol{\delta}(\bar\bu)^\bot ||\leq \bar\ve \leq
d$, where $d$ is the constant from \eref{deltacond}, and let $\xi$
be given by \eref{xibeta}. Then, setting $ \eta_j:= (1-\xi) \bar\ve
2^{-j}$, $\bvv_0:=\bar\bu$, and doing for $j=0,1,2,3,\ldots$
\beqn \label{viterate} \bvv_j - RES(\bvv_j,\eta_j) \to
\widehat{\bvv}_{j+1}, \eeqn one has \beqn \label{deltaeps}
||\boldsymbol{\delta}(\bvv_j)^\bot ||\leq \bar\ve %
,~~~  j=0,1,2,\ldots
\eeqn
and
\beqn
\label{reducedelta}
||\boldsymbol{\delta}(\bvv_j)^\bot ||\leq \xi^j\bar\ve/\beta ,~~~ j=0,1,2,\ldots
\eeqn
\end{lemma}
\noindent
{\bf Proof:}
Setting  $\br_{\eta}:= RES(\bvv_j,\eta)$, we have, by definition,
\begin{eqnarray}
\label{3.1.11}
\widehat\bvv_{j+1}-\bu &=& \bvv_j -\bu -\alpha \br(\bvv_j) + (\alpha\br(\bvv_j)-\br_{\eta_j})\nonumber\\
&=& \bvv_j -\bu -\alpha(\bA(\bvv_j -\bu) -\underline \lambda\bC(\bvv_j-\bu))+\alpha(\lambda(\bvv_j)-\underline \lambda)\bC\bvv_j\nonumber\\
&& ~~~  + (\alpha\br(\bvv_j)-\br_{\eta_j})\nonumber\\
&=& ( \mathbf{I} - \alpha \mathbf{M}  )(\bvv_j-\bu)
+\alpha(\lambda(\bvv_j)-\underline \lambda)\bC\bvv_j +
(\alpha\br(\bvv_j)-\br_{\eta_j}).
\end{eqnarray}
Therefore we obtain, on account of Lemma \ref{Rich(a)}, Lemma \ref{Rich(b)}, and the accuracy of $\br_{\eta_j}$
\begin{eqnarray}
\label{3.1.12}
||\boldsymbol{\delta}(\bvv_{j+1})^\bot ||&\leq& \beta  ||\boldsymbol{\delta}(\bvv_j)^\bot || + \frac{2\alpha ||\bC||\Gamma}{\gamma}\lambda(\bvv_j)||\boldsymbol{\delta}(\bvv_j)^\bot ||^2+ \eta_j\nonumber\\
&=&\Big(\beta  + \frac{2\alpha ||\bC||\Gamma}{\gamma}\lambda(\bvv_j)||\boldsymbol{\delta}(\bvv_j)^\bot ||\Big)||\boldsymbol{\delta}(\bvv_j)^\bot || + \eta_j.
\end{eqnarray}
Now recall from Remark \ref{rem2.1} that $\lambda(\bvv_j)\leq K$ as long as $||\boldsymbol{\delta}(\bvv_j)^\bot ||\leq \sqrt{a}\leq d$ (which is, in particular,
valid by assumption at the initialization step $\bvv_0 =\bar\bu$).

 From
\eref{monotone} we then know that
$$
\|\boldsymbol{\delta}(\bvv_{1})^\bot\|\leq
\Big(\beta  + \frac{\alpha ||\bC||\Gamma}{\gamma}K 
||\boldsymbol{\delta}(\bar\bu)^\bot ||\Big)||\boldsymbol{\delta}(\bar\bu)^\bot || +(1- \xi)\bar\ve
\leq \xi ||\boldsymbol{\delta}(\bvv_0)^\bot || +(1-\xi)\bar\ve \leq \bar\ve,
$$
where $\xi <1$ is given by \eref{xibeta}. Hence,  we conclude that $||\boldsymbol{\delta}(\bvv_1)^\bot ||\leq \bar\ve$. One easily concludes inductively that  \eref{deltaeps} holds (actually with strict inequality for $j\geq 1$).

  More precisely, repeating the reasoning in \eref{3.1.12}, we obtain upon elementary
  calculations and using that $\xi= (1+\beta)/2$,
\beqn
\label{3.1.13}
||\boldsymbol{\delta}(\bvv_{j+1})^\bot ||\leq \xi^{j+1}\bar\ve +\sum_{i=0}^j \xi^{j-i}\eta_i \leq \bar\ve \xi^{j+1}/\beta,
\eeqn
which was to be shown.\feddich

Now recall from Lemma \ref{lem:lema} that, once the orthogonal error component drops below
another possibly smaller threshold
\beqn \label{bardelta} \bar{\delta} := \frac 12\Big\{
\Big(1+\frac{\gamma}{\Gamma K ||\bC|| \,
||\bM_\bot^{-1}||}\Big)^{1/2}-1\Big\},
 \eeqn
the orthogonal error component of the approximate eigenvectors
behaves essentially as the residual.
\begin{remark}\label{rem3.1.2}
If $\bar{\delta} < d$ at most
$$
m_0 = \left\lceil
\frac{\log\Big(\frac{\beta\bar{\delta}}{\overline{\epsilon}}\Big)}{\log
\xi}\right\rceil
$$
iterations of the form \eref{viterate} suffice to provide an approximation
 $\bx_{m_0}$ to $\bu$ satisfying
\beqn \label{initialguess} ||\boldsymbol{\delta}(\bx_{m_0})^\bot ||
\leq \min\,\{d,\bar{\delta}\}. \eeqn
\end{remark}

We are now ready to formulate the main adaptive eigensolver. On account of Remark \ref{rem3.1.2}
we shall assume for simplicity without loss of generality
 that the initial guess already satisfies the somewhat more
stringent accuracy tolerance \eref{initialguess}.\\

\footnotesize
\textit{\textbf{MINIEIG(A,C,$\,\ve$)}}$\to (\lambda(\ve), \bu(\ve))$\hrule
\medskip \noindent
\textit{\textbf{(i) Initialization:}} Choose $\ve_0\leq  \min\,\{d,\bar{\delta}\}$ (see \eref{deltacond}, \eref{bardelta}) and
$|| \mathbf{x}_0 || =1 $ s.t. $\|\boldsymbol{\delta}(\bx_0)^\bot\|\leq \ve_0$;\\
\phantom{ooo}set $\bar\ve := \ve_0$, $\bar\bu := \bx_0$;\\     
  \textit{\textbf{(ii) Iteration block:}}  set $\bar\bu\to \bvv$, $(1-\xi) \bar\ve\to \eta$, $j=0$;\\
  \phantom{ooo}\textit{{do}}\\
 \phantom{oooooo}$ RES(\bvv,\eta)\to \br_\eta$;
  ~~ $ \bvv - \br_\eta \to \widehat\bvv$;
     ~~$\widehat \bvv /||\widehat\bvv|| \to \bvv$;~~$\eta/2\to\eta$, $~ j+1\to j$; \\
 \phantom{ooo}\textit{{until}} $j\geq \log(c_1\beta)/\log \xi$ \textit{or}  $\eta + ||\br_\eta|| \leq c_1\bar\ve/M$ ~\textit{ for a fixed} $~c_1\leq
 \min\,\{2/(5\ve_0), (\sqrt{3}+ 5/2)^{-1}2^{-3/2}\}$;\\
 \textit{\textbf{(iii) Coarsening:}}
  \textit{do }  $\bw = APPROX (\bvv , 3c_1\bar \ve/\sqrt{2})$;\\
 \phantom{ooo}\textit{normalize} $\mathbf{\bar\bu} = \frac{\mathbf{w} }{\|{\mathbf{w} }
\|} $;  \textit{if} $~\bar\ve/2\leq \ve$, \textit{stop and set} $~\lambda(\ve)=RES(\bar\bu,\ve)$,
$~\bu(\ve):= \bar\bu$;\\
\phantom{ooo}\textit{else} $~\bar\ve/2 \to \bar\ve$; ~\textit{go to (ii)} .\\
 \medskip\hrule
\normalsize
\bigskip

Note that a block of perturbed iterations (ii) is interrupted by a coarsening step (iii)
as soon as a threshold criterion is met. This criterion involves actually two alternatives
that will both be seen to be met after a uniformly bounded finite number of steps in (ii).
One of the stopping tests is an {\em a-posteriori} test based on the numerical residual
and the rationale is that it may actually be met earlier than the other test which is based on the
bounds in \eref{reducedelta} which might be too pessimistic, in particular at later stages
when the reduction constants decrease. The role of the coarsening step is to control the
complexity of the overall scheme in a similar way as in adaptive schemes for operator equations
(see e.g. \cite{CDD2,CDD3}).

\begin{theorem}
\label{thconv}
The scheme \textbf{MINIEIG(A,C,$\,\ve$)} terminates for any given $\ve>0$ after finitely many steps
and outputs an approximate eigenpair $ (\lambda(\ve), \bu(\ve))$ with $\bu$ normalized satisfying
\beqn
\label{3.1.10}
||(\bu-\bu(\ve))^\bot|| \leq \ve,~~~~  |\lambda(\ve)- \underline \lambda|\leq \ve,
\eeqn
where $(\underline \lambda, \bu)$ is an exact ground state eigenpair.
\end{theorem}

\begin{remark}
Given $\bu(\ve)$ satisfying the first estimate in \eref{3.1.10} we can approximate $\underline \lambda$
within a tolerance proportional to $\ve^2$ instead of \eref{3.1.10} by
applying $RAYL(\bu(\ve),\ve^2)$ which of course requires a correspondingly higher cost
due to the the required more accurate applications of $SCAL$ and ultimately of $APPLY$, see Section \ref{sect3.3}.
\end{remark}
\noindent
{\bf Proof of Theorem \ref{thconv}:}  To analyze the effect of the various perturbations of the exact iteration suppose that $\bu_k:=\bar\bu$
is the output of the coarsening step after the $k$-th cycle through the coarsening step (iii), i.e.
at this stage we have $\bar\ve =2^{-k}\ve_0$. Moreover, let $\bvv_j$ be the result after $j$ perturbed iterations
in the iteration block (ii) starting with $\bvv_0=\bar\bu=\bu_k$. The corresponding tolerance
$\eta=\eta_j$ is then given by $\eta=(1-\xi)\bar\ve 2^{-j}=(1-\xi)\ve_0 2^{-k-j}$.
Invoking Lemma \ref{lem3.1.2}, ensures that $||\boldsymbol{\delta}(\bvv_j)^\bot ||\leq \bar\ve$
and $||\boldsymbol{\delta}(\bvv_j)^\bot ||\leq \bar\ve \xi^j /\beta$ and hence that the error reduces at least by
a fixed rate. Hence, the stopping criterion in (ii) is met after at most $J:=\left\lceil \log(c_1\beta)/\log \xi\right\rceil$ steps
giving
\beqn
\label{3.1.15}
||\boldsymbol{\delta}(\bvv_J)^\bot || \leq c_1\bar\ve.
\eeqn
 Moreover, the initialization ensures that Lemma \ref{lem:lema} applies
to the iterates $\bvv_j$, which means that $||\boldsymbol{\delta}(\bvv_j)^\bot ||$ is controlled by the residuals.
Thus
\beqn
\label{3.1.14}
||\boldsymbol{\delta}(\bvv_j)^\bot ||\leq M\,||\br(\bvv_j)|| \leq M(||\br_{\eta_j}|| + \eta_j).
\eeqn
On the other hand, also by Lemma \ref{lem:lema},
$$
||\br_{\eta_j}|| ~~\leq~~ \|\br(\bvv_j)\| + \eta_j ~~\leq~~ M' ||\boldsymbol{\delta}(\bvv_j)^\bot || + \eta_j~~\leq~~ M'\bar\ve \xi^j/\beta
+\eta_j.
$$
Thus $||\br_{\eta_j}|| + \eta_j$ will also drop below the threshold $c_1\bar\ve/M$ after finitely many
steps, and in fact, possibly earlier than the alternative step bound $J$, if the bounds in \eref{reducedelta} are overly pessimistic.

In summary, the input $\bvv$  in (iii) at this stage satisfies
$||(\bu -\bvv)^\bot ||\leq c_1\bar\ve$. Assuming without loss of generality that $\lll \bu,\bv\rr > 0$,
we can invoke   Lemma
\ref{projectornorm} to conclude that
\beqn
\label{3.1.16}
|| \bu^\circ -\bvv ||\leq \sqrt{2}c_1\bar\ve  =: \sigma,
\eeqn
where $\bu^\circ := \bu/|| \bu ||$. Since by (iii), we have $\bw =APPROX(\bvv, 3\sigma/2)$, we obtain
by triangle inequality
$
|| \bu^\circ -\bw|| \leq 5\sigma/2$. Moreover, by definition of the routine $APPROX$, we have
$\lll \bw, \bvv -\bw\rr =0$ so that  $||\bw|| \geq (1- (3\sigma/2)^2)^{1/2}$.
Thus, the normalized vector $\bar\bu := \bw/||\bw||$ satisfies with $(1- (3\sigma/2)^2)^{-1/2}=: 1+ g$
\beq
||(\bu -\bar\bu)^\bot ||&\leq &|| \bu^\circ -\bar\bu || \leq ||\bu^\circ -\bw|| + ||\bw - \bar\bu||\\
&\leq &
 \frac{5\sigma}{2} + ||\bw || \Big(\frac{1}{||\bw||} - 1\Big)\leq
 \frac{5\sigma}{2} + ||\bw || \Big(\frac{1}{(1- (3\sigma/2)^2)^{1/2}} - 1\Big)\\
&\leq&
 \frac{5\sigma}{2} + g.
 \eeq
Noting that for $b\leq 1/2$ one has $(1-b^2)^{-1/2}\leq 1+ 2b/\sqrt{3}$, we conclude that $g\leq \sqrt{3}\sigma$, so that
\beqn
\label{3.1.17}
||(\bu -\bar\bu)^\bot ||\leq (\sqrt{3} +5/2)\sigma =\sqrt{2}c_1  (\sqrt{3} +5/2)\bar\ve \leq \bar\ve/2,
\eeqn
by our choice of the constant $c_1$. Thus, in summary we have shown that after a uniformly
bounded finite number of perturbed iterations in (ii) with initial accuracy $\bar\ve$,
one branches into (iii)
whose output is either sufficiently accurate or serves as input for (ii) with improved accuracy $\bar\ve/2$.
Hence the algorithm terminates after finitely many cycles through
(ii), (iii), namely as soon as $2^{-k}\ve_0\leq \ve$.\feddich

\subsection{Complexity estimates}\label{sect3.2}
It remains to analyze the computational complexity of the above scheme {\it \textbf{MINIEIG}}.
The subsequent analysis follows similar lines as used before in connection with
adaptive solvers for operator equations. We formulate an ideal benchmark which describes
the minimal cost needed to achieve a desired accuracy tolerance $\ve$ for an approximate
normalized eigenvector. This lower bound is simply the number of entries needed in any
finitely supported sequence to approximate $\bu$ within accuracy $\ve$ (or equivalently the normalized
$\bu^\circ$ within a fixed factor). This naturally leads to the notion of best $N$-term approximation
that we briefly recall first.

Let $\Sigma_k$ denote the set of all $\bx\in \ell_2(\cI)$ which have at most $k$ nonvanishing entries.
Then
$$
\sigma_N(\bx):= \inf_{\bz\in\Sigma_N}||\bx -\bz||
$$
is the error of best $N$-term approximation in $\ell_2(\cI)$. Then
$$
|\bx|_{\cA^s}:= \sup_{N\in \N} N^s\sigma_N(\bx)
$$
is a (quasi-)seminorm and
$$
\cA^s:= \{\bx\in \ell_2(\cI): ||\bx||_{\cA^s} := ||\bx|| + |\bx|_{\cA^s}< \infty\}
$$
is a (quasi-)Banach space. Thus, for $s>0$, the unit ball of $\cA^s$ is the set of all
those sequences in $\ell_2(\cI)$ whose error of best $N$-term approximation decays at least
as $N^{-s}$ and hence a compact set in $\ell_2(\cI)$. Another way to view this is the following:
In order to approximate a given $\bx\in \cA^s$ with accuracy $\ve$ it takes at most
$N_\ve= \ve^{-1/s}|\bx|_{\cA^s}^{1/s}$ entries to do so, and in the worst case over the whole
unit ball of $\cA^s$ the necessary order of entries is exactly $\ve^{-1/s}$. In what follows this
relation
$$
\mbox{accuracy}~~\ve ~~~~  \leftrightarrow ~~~~  \mbox{degrees of freedom}  ~~
 \ve^{-1/s}|\cdot|_{\cA^s}
$$
reflecting $s$-sparsity of elements in $\ell_2(\cI)$ will be a central orientation
in the subsequent developments.

Let us pause to mention that $\cA^s$ can also be characterized as a {\em weak} $\ell_p$ space.
In fact, denote by $\bx^*$ the nonincreasing rearrangement of $\bx\in \ell_2(\cI)$, i.e.
$x^*_{i+1}=|x_{j_{i+1}}|\leq x^*_i=|x_{j_i}|$, $j=1,2,\ldots$. Then $\ell_p^w(\cI)$
is comprized of all those $\bx\in \ell_2(\cI)$ for which
$$
|\bx|_{\ell_p^w(\cI)}:= \sup_{n\ge1} n^{1/p}x^*_n <\infty.
$$
It is not hard to show that (see \cite{devore98})
$$
\cA^s = \ell_p^w,~~~~  \frac 1p = s + \frac 12.
$$
Moreover, it is easy to see that $\ell_p\subset \ell_p^w$ but for any $q<p$ one has
$\ell_p^w\subset \ell_q$, so that $s$-sparse sequences are almost just $\ell_p$ summable
sequences with $s$ and $p$ related as above.

The first key ingredient of the analysis is the following {\em coarsening lemma} that
explains, in particular, the role of step (iii) in {\it \textbf{MINIEIG}}, see \cite{Cohen,CDD3}.

\begin{lemma}
\label{lem:approx} 
Assume that $ \mathbf{v} \in \cA^s $ for for some $s>0$ and suppose that $\bx \in \ell_2(\mathcal{I})$
is any finitely supported sequence in $\ell_2(\cI)$ such that $ \| \mathbf{v} - \bx
    \|\leq \ve $.
    Moreover fix  $b > 0$ and set
 \beqn
 \label{3.2.1}
 \bw :=\displaystyle{ {\rm argmin}_{||\bv -\bz||\leq (1+ b)\ve }}\#{\rm supp}\,\bz .
 \eeqn

Then there exists a constant $C$ depending only on $b$ and $s$ when $s\to 0$ such that
\beqn
\label{3.2.2}
 \| \bw \|_{\cA^s} \leq
C \| \mathbf{v} \|_{\cA^s},
\eeqn
and
\beqn
\label{3.2.3}
  \# \supp \bw \leq C ||\mathbf{v}||_{\cA^s}^{\frac{1}{s}}\ve^{-\frac{1}{s}},~~ ||\mathbf{v}- \bw \| \leq (2 + b) \ve \leq
\|\mathbf{v}\|_{\cA^s} (\#{\rm supp}\,\bw)^{-s},
\eeqn
 hold uniformly in $
\ve >0$.
\end{lemma}

In other words, thresholding a given finitely supported approximation at a slightly higher tolerance
than the accuracy of approximation provides essentially a best $N$-term approximation to the
(possibly unknown) approximand.

As mentioned before, strictly speaking the cost of determining $\bw$ is essentially \\
$(\#{\rm supp}\,\bx))\log(\#{\rm supp}\,\bx))$. However, at the expense of  a slightly worse
target accuracy than $(1+b)\ve$ one can get away with quasi-sorting based on binary binning
at a computational cost that stays proportional to $(\#{\rm supp}\,\bx))$.
For simplicity, we shall therefore ignore the $\log$-factor in what follows and
use that $APPROX$ can be realized in linear complexity of the input size.

\begin{remark}
\label{rem3.2.1}
Let us denote again by $\bar\bu_k$ the result of the $k$-th application of step (iii) in
\textbf{MINIEIG}. Then, if  the ground state $\bu$ belongs to $\cA^s$ for some $s>0$, 
we have
\beqn
\label{3.2.4}
||\bar\bu_k||_{\cA^s}\lsim ||\bu||_{\cA^s},~~~\#{\rm supp}\,\bar\bu_k \lsim   (\ve_02^{-k})^{-1/s}
||\bu||_{\cA^s}^{1/s},~~~||\bu^\circ  -\bar\bu_k || \leq \ve_0  2^{-k},
\eeqn

uniformly in $ k\in \N$.
\end{remark}

\noindent
{\bf Proof:}  The last relation in \eref{3.2.4} has been already established in the proof of Theorem
\ref{thconv}, see \eref{3.1.17}. The rest is then an immediate consequence of Lemma \ref{lem:approx}.
\feddich

Hence to estimate the computational complexity of \textit{\textbf{MINIEIG}}$(\bA,\bC,\ve)$
depending on $\ve$ as $\ve\to 0$, it remains to bound the computational work in each block (ii).
\begin{proposition}
\label{prop3.2.2}
Assume that for some $s>0$ one has $\bu \in \cA^s$, and that $\br_\eta := RES(\bvv,\eta)$ satisfies
\newcommand{\bound}{|| \bu ||_{{\cA}^s} + || \bvv ||_{{\cA}^s} } 
\begin{equation}
\label{3.2.5}
\begin{array}{rcl}
|| \br_\eta ||_{{\cA}^s} & \leq & C( \bound ), \\[2mm]
 \#{\rm supp}\,\br_\eta, ~ \#{\rm flops}(\br_\eta)  & \leq & C
\eta^{-1/s}\big(||\bu||_{\cA^s}^{1/s} +||\bvv||_{\cA^s}^{1/s}\big),
\end{array}\end{equation}

holds for some constant $C$ independent of $\eta$.  Here $\#{\rm flops}(\br_\eta)$ denotes the number
of arithmetic operations needed to compute $\br_\eta$.
Then
the output $(\lambda(\ve),\bu(\ve))$ of \textbf{MINIEIG}$(\bA,\bC,\ve)$ satisfies in
addition to \eref{3.1.10}
\beqn
\label{3.2.6}
 \#{\rm supp}\,\bu(\ve), \, \#{\rm flops}(\bu(\ve)) \leq C
\ve^{-1/s}||\bu||_{\cA^s}^{1/s},~~~||\bu(\ve)||_{\cA^s} \leq C||\bu||_{\cA^s},
\eeqn

for some constant $C$ independent of $\bu$ and $\ve$.
\end{proposition}

\noindent
{\bf Proof:} Applying if necessary the coarsening lemma to the initial guess
(starting from a correspondingly slightly higher initial  accuracy) we have
$|\bx_0|_{\cA^s}\lsim |\bu|_{\cA^s}$, $ \#{\rm supp}\,\bx_0 \lsim \ve_0^{-1/s}$.
 Suppose now that at the $k$-th stage of (ii)
the input $\bar\bu =\bar\bu_{k-1}$ satisfies for some constant $C$ independent of $k$
\beqn
\label{3.2.7}
\#{\rm supp}\,\bar\bu, \, \#{\rm flops}(\bar\bu) \leq C
\bar\ve^{-1/s}(||\bu||_{\cA^s}^{1/s} + ||\bar\bu||_{\cA^s}^{1/s}).
\eeqn
By assumption \eref{3.2.5}, the same estimates hold for the intermediate iterates $\bvv_j$
in (ii) (with $\bvv_0=\bar\bu_{k-1}$) with constants however now depending on $j$.
As shown in the proof of Theorem \ref{thconv} each block (ii) has at most a uniformly bounded
number $J$ of iterations, so that \eref{3.2.7} still holds for the input to step (iii) with some constant
$C=C(J)$, $J$ being the upper bound of the number of iterations in (ii).
As pointed out in Remark \ref{rem3.2.1}, the application of $APPROX$ produces then
$\bar\bu=\bar\bu_k$ satisfying again \eref{3.2.7} for $\bar\ve_k=\bar\ve_{k-1}/2$
with a constant coming from the coarsening
lemma that is independent of $k$. Thus summing the computational cost of each block (ii)
and using a straightforward geometric series argument confirms the claim. \feddich

Thus, it remains to verify \eref{3.2.5} for the approximate residuals.
According to the ingredients of the approximate residuals provided by the routine $RES$,
given in the previous section, the key issue here is the efficient approximate application of
the matrices $\bA$ and $\bC$ through the routine $APPLY$. It is well known by now that for
a wide range of (local and global) operators and suitable wavelet bases the corresponding
operator representations in wavelet coordinates are {\em nearly sparse}.
A precise formulation of this property   that is fulfilled
in many concrete cases, reads as follows, see e.g. \cite{Denc,DHS}.\\

\noindent
{\bf $s^*$-compressibility:} {\it  Let $s^*$ be a positive real. $\cB^{s^*}$ denotes the set of matrices (over
$\cI\times \cI$) with the following properties: $\bB\in \cB^{s^*}$ if for every $k \in \N$, there is a matrix $\mathbf{B}_k$ with at most $2^k$ entries in each row and column satisfying $||\mathbf{B} - \mathbf{B}_k|| \leq C2^{-ks^*}\alpha_k,$ where $(\alpha_k)_{k=0}^{\infty}$ is an $\ell_1$-sequence of positive numbers.
Elements in $\cB^{s^*}$ are called $s^*${\em -compressible}.}

One can show that any $\mathbf{B}\in \cB^{s^*}$ maps the $\cA^s$-spaces  boundedly into themselves for $ s < s^*$, cf. \cite{CDD1}, Section 3. For elements in $\bB\in \cB^{s^*}$ a concrete realization
of $APPLY(\bB,\cdot,\cdot)$ has been
developed in \cite{CDD1} whose properties are given in the following lemma.

\begin{lemma}
\label{dahmenlem}
Assume that $\bB\in \cB^{s^*}$.
Given a tolerance $\delta > 0$ and a vector $\bx$ with finite support, the algorithm $APPLY (\mathbf{B}, \bx, \delta)$ produces a vector $\mathbf{w}$ which satisfies
$$
||\mathbf{B}\bx - \mathbf{w}|| \leq \delta.
$$

Moreover, for $s<s^*$ one has:
 \begin{enumerate}
\item[(i)] The output vector $\mathbf{w}$ satisfies
$$
||\mathbf{w}||_{\cA^s} \lesssim ||\bx||_{\cA^s}; \qquad \supp \mathbf{w} \lesssim
||\bx||_{\cA^s}^{\frac{1}{s}} \delta^{-\frac{1}{s}}.
$$
\item[(ii)]
The number of entries of $\mathbf{B}$ to be computed to obtain $\bw$
 is $\lesssim ||\bx||_{\cA^s}^{\frac{1}{s}} \delta^{-\frac{1}{s}};$.
\end{enumerate}
\end{lemma}
For the \textit{proof}, again see \cite{CDD1}.

If for $s<s^*$ the number of arithmetic operations needed to compute $\bw$
does not exceed $C||\bx||_{\cA^s}^{\frac{1}{s}} \delta^{-\frac{1}{s}} + 2\supp \bx$ (using quasi-sorting instead of exact sorting which would entail an additional $\log$-factor),
the matrix $\bB$ is called $s^*$-{\em computable}, \cite{stev1,stev2}. For the verification of
$s^*$-computibility for a wide class of  operators, see \cite{stev1,stev2}.\\

\noindent
{\bf $s^*$-sparsity:} {\it The matrix $\bB$ is called {\em $s^*$-sparse} if there exists a scheme
$APPLY$ satisfying the properties listed in Lemma \ref{dahmenlem} whose computational complexity
remains proportional to the output size.}\\

Clearly $s^*$-computable matrices are examples of $s^*$-sparse matrices. It is important
to note though that there are further important examples. Recalling the form of $\cL$ in remark \ref{remC},
the matrix $\bA$ may actually be the product of several matrices. An $APPLY$-scheme for such products can easily be obtained by composing individual $APPLY$-schemes designed e.g.
for compressible matrices, see \cite{CDD2} for the treatment of least squares formulations.

Another important case concerns the application of $\bC$ in a fairly
general setting.

Actually $ \bC $  represent the {\em inverse} $\bC=\mathbf{S}^{-1}$
of an (elliptic) operator $\cS$ in  the sense of (\ref{Lelliptic}).
Thus the $APPLY$-scheme for $\bC$ may just mean the adaptive
approximate solution of an operator equation for an $ \cH$-elliptic
operator $\cS: \cH \to \cH'$. In many cases its output has been
shown to satisfy the properties in Lemma \ref{dahmenlem} and thus
gives rise to an $s^*$-sparse $APPLY$-scheme.

In summary, it is important to keep in mind
that the actual realizations of $APPLY$ for $\bA$ and $\bC$ may be completely different
but should satisfy the properties of Lemma \ref{dahmenlem}.

\begin{assumption}
 \label{sparsityass}
 The matrices
$\bA $ and $\bC$ are $s^*$-sparse for some $s^*>0$ (For $\bC$ this is trivially the case when
$\bC$ can be applied exactly in linear time).
\end{assumption}

The main result of this paper can now be formulated as follows.
\begin{theorem}
\label{thmain}
Assume that $\bA$ and $\bC$ are $s^*$-sparse for some $s^*>0$ and that
the parameters $\alpha, \beta$ are chosen according to Lemma \ref{Rich(a)}.
Then for any $\ve >0$, the scheme \textbf{MINIEIG}$(\bA,\bC,\ve)$ after
finitely many steps produces an approximate eigenpair $(\lambda(\ve), \bu(\ve))$
with normalized $\bu$ satisfying
\beqn
\label{3.2.8}
||(\bu-\bu(\ve))^\bot||\leq \ve,~~~~  |\underline \lambda-\lambda(\ve)|\leq \ve,
\eeqn
where $(\underline \lambda,\bu)$ is the exact ground state solution of \eref{evp}.

Moreover, if $\bu\in \cA^s$ for some $s<s^*$, then one has
\beqn
\label{3.2.9}
\#{\rm flops}\,\bu(\ve),\, \#\supp\,\bu(\ve)\lsim \ve^{-1/s}||\bu||_{\cA^s}^{1/s},~~~~
||\bu(\ve)||_{\cA^s}\lsim ||\bu||_{\cA^s},
\eeqn
where the constants are independent of $\ve$ and $\bu$ but depend only on $s$ when $s$ approaches
$s^*$.
\end{theorem}
\noindent
{\bf Proof:} \eref{3.2.8} has been already shown in Theorem \ref{thconv}.
To prove the rest of the claim we employ Proposition \ref{prop3.2.2} which requires
confirming the property \eref{3.2.5}.  Towards this end, recall from Section
\ref{adaptalgo} that the realization of $RES(\bx,\eta)$ requires the approximate application of
$\bA,\, \bC$ within tolerances that are uniformly  bounded from below and above by fixed multiples of $\eta$ (see \eref{3.1.9}) as well as the computation of $RAYL(\bx,\eta/(4\alpha))$.
Now, this latter routine requires the evaluation of the routine $SCAL$ for $\bA$
and, unless $\bC$ can be applied exactly, for $\bC$. By Remark \ref{rem3.1.1} and Lemma \ref{lemrayl},
the tolerances needed in the $SCAL$ routines remain uniformly proportional to the target accuracy in
$RAYL(\bx,\eta/(4\alpha))$ which is proportional to $\eta$. Since all these routines rely on the scheme $APPLY$
with respect to tolerances proportional to $\eta$, the relations \eref{3.2.5} follow from
Lemma \ref{dahmenlem}. This finishes the proof. \feddich

\subsection{Possible quantitative improvements}\label{sect3.3}
The scheme analyzed above should be viewed as one possible realization of an adaptive
strategy.
To achieve quantitative improvements of schemes of the above type one might try to exploit the fact
that the (exact) Rayleigh quotient exhibits essentially the square accuracy of the corresponding
approximate eigendirection. We shall only sketch such  strategies whose details would essentially follow
 analogous lines as discussed above.

Recall that the evaluation of the scalar products requires the tightest accuracy tolerances
(although still proportional to the target accuracy of $RES$). Thus there are two angles that might help
to reduce computational complexity, namely:\\

a) Trying to speed up the calculation of scalar products;

b) Reducing the number of calls of $RAYL$.

\subsubsection{Fast evaluation of scalar products}

As for a), the naive approach outlined above rests on the accurate calculation of a complete
matrix/vector product although one computes at the end only a single number. We have already
mentioned a possible approach in \eref{3.1.3} that might help reducing the computational complexity.

Recalling that an approximate eigendirection of accuracy $\ve$ would
give rise to a (precise) Rayleigh quotient that approximates
$\underline{\lambda}$ within accuracy of the order of $\ve^2$ (see
Lemma \ref{Rich(b)}), a first natural idea is to postprocess the
result of {\it \textbf{MINIEIG}}$(\bA,\bC,\ve)$, satisfying
\eref{3.2.8},  so as to obtain an approximation to $\lambda$ of
order $\ve^2$. In fact, in view of Lemma \ref{Rich(b)},  it would
make sense to compute
$$
\lambda^*(\ve)= RAYL(\bu(\ve),\ve^2),\quad \mbox{so that} ~~
|\lambda^*(\ve)-\underline\lambda|\leq (1+C)\ve^2, ~~\mbox{with}~~  C:= 2\Gamma K/\gamma,
$$
see \eref{ewzeile2} and \eref{raylabs}. However, using the simple
version of the routine $SCAL$ described below \eref{3.1.2}, the
computational complexity could be of the order of $\ve^{-2/s}$ when
$\bu \in \cA^s$. Let us now point out that this cost can be reduced
significantly by employing more refined versions of $SCAL$ along the
lines of \eref{3.1.3}. To this end, suppose that $\bB$ is an
$s^*$-compressible matrix and recall e.g. from \cite{CDD1,DHS} the
following key idea of constructing an $APPLY$ scheme for an
approximate matrix/vector computation. Fix any $\bar s <s^*$. Given
$\zeta >0$ and $\bvv$ of finite support, let $\bvv_\zeta :=
APPROX(\bvv, \zeta)$ and set
$$
\bvv_{-1}(\zeta) := \bvv -\bvv_\zeta, \quad \bvv_j(\zeta):=
\bvv_{2^{j\bar s}\zeta} - \bvv_{2^{(j+1)\bar s}\zeta}, \quad  j=0
\ldots , J(\zeta):=\lceil \log_2 (||\bvv||/\zeta)/\bar s \rceil \, .
$$
Then, one has
\beqn
\label{3.3.2}
\sum_{j=-1}^{J(\zeta)} \bvv_j(\zeta)
=\bvv, \quad ||\bvv_{-1}(\zeta)||\leq \zeta, \quad  || \bvv_j(\zeta)|| \leq (1+2^{\bar s})2^{j\bar s}\zeta,
\,\, j=1,\ldots, J(\zeta) .
\eeqn
Moreover, it is shown in \cite{DHS} that (when $\bB_j$ are the compressed versions of
$\bB\in \cB^{s^*}$ from the definition of
$s^*$-compressibility)
\beqn
\label{3.3.1}
\bw_\zeta:= \sum_{j=0}^{J(\zeta)}\bB_j \bvv_j(\zeta)\quad \mbox{satisfies}\quad ||\bB\bvv -\bw_\zeta
||\leq C' \zeta,
\eeqn
for some uniform constant $C'$ independent of $\zeta$. For simplicity we shall work with $C'=1$
which can always be arranged through the definition of the $\bB_j$ or by replacing $\zeta$
by $c\zeta$ in the decomposition of $\bvv$.

The announced improved version of $SCAL$ is based on the following observation.
\begin{proposition}
\label{prop3.3.1}
Assume that $\bB\in \cB^{s^*}$ and $\bar s< s^*$ is fixed.
Setting
\beqn
\label{3.3.4}
\epsilon_j :=   \alpha_j (1 + 2^{\bar s})^{-1} 2^{-\bar s j} \delta^{1/2},\quad j=0,\ldots,J(\sqrt{\delta}),\,\,
\epsilon_{-1}:= \alpha_{-1}\sqrt{\delta},\,\, \sum_{j=-1}^\infty \alpha_j =1,
\eeqn
where the summable (to one)  coefficients $\alpha_j$ are again from the definition of $s^*$-compressibility and have algebraic decay, and defining
\beqn
\label{3.3.5}
s(\delta,\bvv):= \sum_{j=-1}^{J(\sqrt{\delta})} \lll \bvv_{j}(\sqrt{\delta}), \bw_{\epsilon_j}\rr ,
\eeqn
where the $\bw_{\epsilon_j}$ are given by \eref{3.3.1}, one has
\beqn
\label{3.3.6}
|s(\delta,\bvv) -\lll \bB\bvv,\bvv\rr|\leq \delta.
\eeqn

Moreover,  whenever ${\rm supp}\,\bvv \leq C \delta^{-1/2s}|\bvv|_{\cA^s}^{1/s}$ for some $s\leq \bar s<s^*$, then
the number of operations $ops(\bvv,\delta)$ needed to compute $s(\delta,\bvv)$
is bounded by a constant multiple of $\Big( |\bvv|_{\cA^s}\delta^{-\frac{1}{2s}}\Big) ^{ \frac{\bar s +2s}{s+\bar s}}$, where the constant
depends only on $C$.
\end{proposition}
\noindent
{\bf Proof:}
By \eref{3.3.2} we can write
$$
|s(\delta,\bvv)-\lll\bB\bvv,\bvv\rr| =\Big| \sum_{j=-1}^{J(\sqrt{\delta})}\lll  \bvv_j(\sqrt{\delta}), \bw_{\epsilon_j} -\bB\bvv\rr
\Big|\leq \sum_{j=-1}^{J(\sqrt{\delta})}\alpha_j\delta\leq \delta ,
$$
where we have used Cauchy-Schwarz and \eref{3.3.4} in the last step, confirming \eref{3.3.6}.

As for the work count, we shall estimate now first the number of operations needed to compute a single
entry of $\bw_\zeta$.
Now   note that (since $\bvv\in \cA^s$ for any $s>0$) for $\Delta_j := {\rm supp}\, \bvv_j(\zeta)$ one has
\beqn
\label{3.3.3}
 \# \Delta_j \lsim (2^{j\bar s}{ \zeta})^{-1/s} |\bvv|_{\cA^s}^{1/s}.
\eeqn
Further, recall that each row of $\bB_j$ has at most $2^j$ entries, so that the computation of one entry of a contribution $\bB_j \bvv_j({\zeta})$ takes, in view of \eref{3.3.3}, at most $2^j$ operations as long as
$j\leq j^*$ when $j^*=j^*(\zeta)$ is the largest integer for which
\beqn \label{3.3.7} 2^{j^*}\leq |\bvv|_{\cA^s}^{1/s} 2^{-j^*\bar
s/s} \zeta^{-1/s}\quad \Longleftrightarrow \quad j^*=\left\lfloor
(s+\bar s)^{-1}\log_2 \Big(\frac{|\bvv|_{\cA^s}}{\zeta}\Big)
\right\rfloor. \eeqn
Thus the computation of a single entry of the partial sum $\sum_{j=0}^{j^*}\bB_j\bvv_j(\zeta)$
takes the order of
\beqn
\label{3.3.8}
2^{j^*}\lsim \zeta^{-1/(s+\bar s)} |\bvv|_{\cA^s}^{1/(s+\bar s)}.
\eeqn
Likewise the number of computations required for a single entry of the remaining sum
$\sum_{j=j^*(\zeta)}^{J(\zeta)} \bB_j\bvv_j(\zeta)$ is, by \eref{3.3.3}, of the order of
\beqn
\label{3.3.9}
\zeta^{-1/s} |\bvv|_{\cA^s}^{1/s} \sum_{j=j^*}^J 2^{-j\bar s/s}
\lsim \zeta^{-1/s} |\bvv|_{\cA^s}^{1/s} 2^{-j^*\bar s/s}\lsim |\bvv|_{\cA^s}^{1/(s+\bar s)} \zeta^{-1/(s+\bar s)}.
\eeqn
Therefore
\beqn
\label{3.3.10}
\# ops(\mbox{for computing one entry of ~}
\bw_\zeta)\lsim \zeta^{-1/(s+\bar s)} |\bvv|_{\cA^s}^{1/(s+\bar s)}.
\eeqn
We shall now estimate the work required by the computation of $\lll \bvv_{(j)},\bw_{\epsilon_j}\rr$.
Note first that
$$
\#{\rm supp}\, \bvv_{-1}(\sqrt{\delta})\lsim \delta^{-1/2s}|\bvv|_{\cA^s}^{1/s},\quad
\#{\rm supp}\, \bvv_{j}(\sqrt{\delta})\lsim 2^{-\bar s j/s}\delta^{-1/2s}|\bvv|_{\cA^s}^{1/s}.
$$
Hence, using \eref{3.3.10} with $\zeta =\epsilon_j$,
  the computation of $\lll \bvv_{j}(\sqrt{\delta}),\bw_{\epsilon_j}\rr$
in \eref{3.3.5}
takes, in view of \eref{3.3.5},  the order of
\beqn
\label{3.3.11}
\#{\rm supp}\, \bvv_{j}(\sqrt{\delta}) |\bvv|_{\cA^s}^{1/(s+\bar s)} \epsilon_j^{-1/(s+\bar s)}~\leq~ \alpha_j^{-\frac{1}{s+\bar s}}
 (1+2^{\bar s})^{\frac{1}{s+\bar s}}|\bvv|_{\cA^s}^{\frac{2s+\bar s}{s(s+\bar s)}}2^{-j\bar s^2/(s\bar s+\bar s^2)}\delta^{-\frac{1}{2s}\Big( \frac{\bar s +2s}{s+\bar s}\Big)}
\eeqn
operations. Summing over $j$ and recalling that the $\alpha_j$ decay polynomially completes the
proof. \feddich

\begin{corollary}
\label{cor3.3.1} Assume that the hypotheses of Theorem \ref{thmain}
are valid. Given $\ve >0$, let $\bu(\ve)$ be the output of
\textbf{MINIEIG}$(\bA,\bC,\ve)$ and let $\lambda^*(\ve)=
RAYL(\bu(\ve),\ve^2)$ where $RAYL$ is based on a version of $SCAL$
derived in an obvious manner from \eref{3.3.5}.  Then one has \beqn
\label{3.3.12} |\underline{\lambda}  
 -\lambda^*(\ve)| \leq (1+2\Gamma K/\gamma))\ve^2,
\eeqn
and, whenever $\bu \in \cA^s$ for some
$s\leq \bar s$, the computational complexity $\# ops(\lambda^*(\ve))$ of $\lambda^*(\ve)$ remains bounded by
\beqn
\label{3.3.13}
\# ops(\lambda^*(\ve))\lsim  |\bu|_{\cA^s}^{\frac{2s+\bar s}{s(s+\bar s)}}\ve^{-\frac{1}{s}\Big( \frac{\bar s +2s}{s+\bar s}\Big)},
\eeqn
where the constant is independent of $\bu$ and $s$.
\end{corollary}
{\bf Proof:} The estimate \eref{3.3.12} is an immediate consequence of \eref{ewzeile2}
and the first relation in \eref{3.2.8}. The complexity estimate, in turn, follows from
Proposition \ref{prop3.3.1} applied to $\bvv =\bu(\ve)$ with $\delta := \ve^2$ together with \eref{3.2.9}.
\feddich

Since $g(s):=(\bar s +2s)/(\bar s +s)$ increases in $s$ and $g(\bar s) = 3/2$, the computational complexity of computing $\lambda^*(\ve)$, and hence the smallest eigenvalue $\underline \lambda$ within
a tolerance of order $\ve^2$, grows at most like $\ve^{-3/2s}$, which is of course much better
than the cost $\ve^{-2/s}$ that would result from applying the original
simple version of $SCAL$. In fact, when $s$ is very small one almost recovers cost of $\ve^{-1/s}$
needed to approximate $\bu\in \cA^s$ within tolerance $\ve$.

Note also that the coarsened versions of $\bvv$ needed in the computation of
$\bw_{\epsilon_j}$ are essentially the same as those in \eref{3.3.5}, so that they can
be reused.
Nevertheless, this version of $RAYL$ is quantitatively still more involved
as the original simple one based on approximating $\bB\bvv$ at the
respective target accuracy. Since these vectors are needed anyway
in the course of {\it \textbf{MINIEIG}} it seems preferable to use
the latter more elaborate version only at the end once $\bu(\ve)$ has
been obtained.

\subsubsection{Modified iterations}
As for b), another option is to modify the ideal iteration itself, again trying to exploit
the fact that the approximation rate of the (exact) Rayleigh quotients is
faster than that of the eigendirections. We only give a very rough sketch of
the idea. Instead of the final accuracy one could run {\it \textbf{MINIEIG}} first with
some intermediate accuracy $\ve$ outputting again $\bu(\ve)$ satisfying \eref{3.2.8}.
One could then fix $\bar\lambda =\lambda^*(\ve)$, set $\bx_0:= \bu(\ve)$
and iterate
$$
\bx_{i+1}=\bx_i-\alpha(\bA -\bar\lambda\bC)\bx_i,
$$
so that, by straightforward calculations, one obtains
$$
(\bx_{i+1}-\bu)^\bot  = (\bx_i -\bu)^\bot  -\Phi(\bx_i -\bu)^\bot + \alpha (\bar\lambda -\underline \lambda)
(\mathbf{I-P)}\bC \bx_i,
$$
and hence
by recursion
\beqn
\label{3.3.14}
|| (\bx_i-\bu)^\bot || \leq \beta^i\ve + \frac{(1+2K\Gamma\alpha/\gamma)\ve^2}{1-\beta}||\bC||.
\eeqn
Thus, in principle, after $|\log \ve |$ steps one has quadratic accuracy of $(\bx_i-\bu)^\bot$
without further applications of $RAYL$.
However, the iterates $\bx_i$ are no longer orthogonal to the residuals
$\bar\br_i := (\bA-\bar\lambda \bC)\bx_i$. Observing that
$\lll \bar\br_i,\bx_i\rr = (\lambda(\bx_i)-\bar\lambda)\lll\bC \bx_i,\bx_i\rr$, we have
$$
||\bx_{i+1}||^2 = ||\bx_i ||^2+\alpha^2 ||\bar\br_i ||^2 -2\alpha
(\lambda(\bx_i)-\bar\lambda) ||\bx_i||_\bC^2  \geq ||\bx_i
||^2+\alpha^2 ||\bar\br_i ||^2 -2\alpha
(\lambda(\bx_i)-\bar\lambda)||\bC||  ||\bx_i||^2,
$$
so that we cannot guarantee any more that the iterates have
increasing norms. Nevertheless, since $|\lambda(\bx_j)-\bar\lambda|$
is expected to be of the order $\ve^2$ and $||\bx_0||=1$, there is
still a fixed positive constant $b\geq 1/2$ say, so that after
$J\sim |\log \ve|$ steps, one still has $||\bx_J|| \geq b$, so that
renormalizing the $J$th iterate preserves quadratic accuracy. This
way, one expects to catch up quadratic accuracy of the approximate
eigendirections without any intermediate computation of Rayleigh
quotients. One can then continue such a block of iterations with
initial guess $\bx_J/||\bx_J||$ that approximates $\bu$ now within a
tolerance of the order of $\ve^2$ so that a corresponding
approximate eigenvalue can be computed within tolerance $\ve^4$ with
the aid of the above
fast computation of scalar products, etc. We leave the details to the reader.\\

Concerning optimal line search \eref{linesearch} and subspace acceleration
techniques, we mention that this requires the solution of a two-dimensional
(or more generally small)
generalized eigenvalue problem. However, setting up the corresponding matrices
requires the computation of scalar products. Since, generally, we cannot compute these quantities exactly, a careful assessment of the perturbation effects
would be necessary. Also, the possibilities for realizing enhanced accuracy at reduced cost mentioned above might be relevant in this context.


\section{Remarks concerning an application -\\ The Schrödinger equation}
\label{schroedinger}

We conclude with some brief comments concerning a scenario that has actually
motivated part of this work.\\

We consider the time-independent electronic Schrödinger equation
$H\psi = E_0 \psi$, where $H: H^1( (\R^{3} \times \{\pm
\frac{1}{2}\})^N ) \rightarrow H^{-1}( (\R^{3} \times \{\pm
\frac{1}{2}\})^N )$ is the Hamilton operator of the molecular
$N$-electron system under consideration, and $E_0$ denotes the
lowest eigenvalue of $H$ corresponding to the ground state of the
system. The spectrum of $H$ is real and bounded from below by a
certain $\mu \in \R$ (cf. e.g. \cite{yser} and references therein),
so that the shifted Hamiltonian $H':= H-\mu\cE$ satisfies the
conditions \eref{Lsymmetric} and \eref{Lelliptic} of Section
\ref{sect1}. Letting $\Phi =\{\phi_{\nu}|\nu \in \mathcal{I}\}$ be
an orthonormal basis of $L_2( (\R^{3} \times \{\pm \frac{1}{2}\})^N
)$, we have to solve the eigenvalue equation \beqn {\mathbf H}
\mathbf{c} = E_0{\mathbf c} \qquad \label{schro} \eeqn where the
entries $\mathbf H_{\nu', \nu}$ of the ``matrix Hamiltonian'' in this
formulation are given by
$$
{\mathbf H_{\nu', \nu}} =\langle H\psi_{\nu'}, \psi_{\nu} \rangle.
$$ In view of our discussion in Section \ref{sect1},
the problem is posed in $\ell_2(\cI)$, which gives the so-called
\emph{complete CI} formulation for the Schrödinger equation. To
obtain a Riesz basis for $H ^1$, we can utilize a reference
operator $\mathcal{F}= \sum_{i=1}^N F_i $, where the $ F_i : H^1
(\mathbb{R}^3 \times \{\pm \frac{1}{2}\} ) \to  H^{-1} (
\mathbb{R}^3 \times \{\pm \frac{1}{2}\} ) $ are single particle
operators with a complete eigenvalue basis. The self consistent Fock
operator from the Hartree Fock equations, see e.g. \cite{helg}, can
be modified to fit  the present purpose. This operator
  is  also bounded as an
operator $H^1  ( ( \mathbb{R}^3 \times \{\pm \frac{1}{2}\}^N)  \to
H^{-1} (( \mathbb{R}^3 \times \{\pm \frac{1}{2}\})^N) $ and its
spectrum is bounded from below,
  so that $ \mathcal{F}' := \mathcal{F} - \mu \mathcal{E} $
  satisfies (\ref{Lelliptic}).
In this case the eigenfunctions $\phi_{\nu}$ and the eigenvalues $\sigma_{\nu}$ of $ \mathcal{F},
\mathcal{F}' $ can easily be computed from those of the single particle
operator $F_i$,
$$\mathcal{F}' \phi_{\nu} = \sigma_{\nu} \phi_{\nu} \ , \ \nu \in
\mathcal{J}.$$
If we choose
  $$\Psi
=  \{\psi_{\nu} := \sigma_{\nu}^{-\frac{1}{2}} \phi_{\nu} |\nu \in
\mathcal{I}\}$$
 this basis satisfies \eref{ne}.
The generalized eigenvalue problem $(\ref{evp})$ then resembles the
symmetry-transformed variant of the Schrödinger equation with \beqnr
  \bC: = \big(\langle \sigma_{\nu}^{-\frac 12}\phi_{\nu},
  \sigma_{\nu'}^{-\frac 12}\phi_{\nu'} \rangle \big)_{\nu ,\nu' \in \cI} \quad
 \makebox{ and }\quad \mathbf{A}:=
 \mathbf{C}^{\frac{1}{2}}\mathbf{H'}\mathbf{C}^{\frac{1}{2}} =
 \mathbf{C}^{\frac{1}{2}}(\mathbf{H} - \mu \mathbf{I})
 \mathbf{C}^{\frac{1}{2}},\label{matrixbezeichnungen}
 \eeqnr
where the conditions \eref{Aund2} and \eref{raylminl2} are valid.

With the above notation, $||\mathbf{x}||_{\mathbf C}$ then gives the
$L_2( (\R^{3} \times \{\pm \frac{1}{2}\})^N )$-error of the wave
function $\varphi:= \sum_{\nu \in \cI} x_{\nu} \psi_{\nu} =
\sum_{\nu \in \cI} c_{\nu} \phi_{\nu}$. For estimating the $H^1(
(\R^{3} \times \{\pm \frac{1}{2}\})^N )$-error of the wave function
$\varphi$, note that the canonical $H^1( (\R^{3} \times \{\pm
\frac{1}{2}\})^N ) $-norm is equivalent to the norms induced by the
inner products $\langle \langle .,. \rangle \rangle := \langle H'
\cdot,\cdot \rangle_{L_2}$ and $\langle \mathcal{F}'\cdot
,\cdot\rangle_{L_2}$, defined by the Hamiltonian and the reference
operator respectively, which gives \beq
 ||\mathbf{x}||_{\mathbf{A}}^2 &=&   \langle \mathbf{C}^{\frac{1}{2}}\mathbf{H'~C}^{\frac{1}{2}}
 \mathbf{x},\mathbf{x}\rangle   \qquad~~=
 ~~\langle \mathbf{H'~C}^{\frac{1}{2}}\mathbf{x},\mathbf{C}^{\frac{1}{2}}\mathbf{x}\rangle \\
 &=& \sum_{\nu \in \cI} \sum_{\nu' \in \cI} c_{\nu} c_{\nu'} \langle H' \phi_{\nu}, \phi_{\nu'} \rangle
  ~=~~ \langle H'~\varphi, \varphi \rangle ~\simeq~ ||\varphi||_{H^1}^2,
 \eeq
and
 \beq ||\mathbf{x}||_{2} ~=~
 ||\mathbf{C}^{\frac{1}{2}}\mathbf{x}||_{\mathbf{C}^{-1}}  = ||\mathbf{c}||_ {\mathbf{C}^{-1}}
 = ||\varphi||_{\mathcal{F}'} \simeq ||\varphi||_{H^1}.
\eeq This  serves also as a typical example where $\mathbf C$ is not
coercive on $\ell_2(\cI)$ (because $||.||_\mathbf C$ is equivalent
to the $L_2$-norm), cf. Remark \ref{remC}.
 The norms $||.||_{\mathbf{A}}$ and $||.||_2 \simeq ||.||_{\mathbf{A}}$ on $\ell_2(\cI)$
 can now be used to estimate the convergence of the $CI$ solution with respect to the $H^1$-norm.


\end{document}